\documentclass[11pt]{amsart}

\newtheorem{theorem}{Theorem}[section]
\newtheorem{proposition}[theorem]{Proposition}
\newtheorem{lemma}[theorem]{Lemma}
\newtheorem{corollary}[theorem]{Corollary}

\theoremstyle{definition}
\newtheorem{remark}[theorem]{Remark}
\newtheorem{definition}[theorem]{Definition}

\usepackage{graphics}
\usepackage{t1enc}
\usepackage{textcomp}
\usepackage{amsmath, amstext,amscd, amsthm}
\usepackage{epic} \usepackage{eepic}
\usepackage{amssymb}

\def\Spec{{\rm Spec}\,}
\def\P{\mathbb{P}}
\def\tens{\otimes}
\def\O{\mathcal{O}}
\def\Hom{{\rm Hom}}
\def\End{{\rm End}}
\def\Aut{{\rm Aut}}

\def\Ext{{\rm Ext}}
\def\wX{\widetilde{X}}
\def\wH{\widetilde{\mathcal{H}}}

\makeatletter

\@addtoreset{equation}{section} \makeatother

\begin{document}

\baselineskip=15.5pt

\title[Frobenius inverse image of the semi-stable boundary]{Frobenius inverse image of the semi-stable boundary in the moduli space of vector bundles}

\author[L. Ducrohet]{Laurent Ducrohet}

\address{CMLS, Ecole Polytechnique, 91128 Palaiseau Cedex, France}

\email{ducrohet@math.polytechnique.fr}

\date{}

\maketitle

\section{Introduction}

Let $X$ be a genus $g\geq 2$ proper and smooth curve over an algebraically closed field $k$ of characteristic $p\neq 2$. We let $J_X$ denote the Jacobian of $X$. We choose once for all a theta characteristic $\kappa_0$ over $X$ and we let $\Theta$ denote the corresponding symmetric principal divisor on $J_X$ representing the canonical polarization of $J_X$. The coarse moduli space $S_X$ of (S-equivalence classes of) rank 2 and trivial determinant semi-stable bundles over $X$ identifies with
a normal subvariety of dimension $4(g-1)-1$ of the linear space $|2\Theta|\cong \mathbb{P}^{2^g-1}$ of effective divisors linearly equivalent to $2\Theta$ on $J_X$. Furthermore, its semi-stable boundary identifies with the Kummer variety $K_X$, which is canonically embedded in $|2\Theta|$.

Assume that $k$ has positive characteristic $p\geq 3$ and denote by $X_1$ (resp. $J_{X_1}$, $S_{X_1}$,...) the Frobenius twist of $X$ (resp. $J_{X}$, $S_{X}$,...). The relative Frobenius $F : X \to X_1$ is a $k$-morphism that induces by
 pull-back a morphism $V_J : J_{X_1} \to J_X$ classically called Verschiebung, étale when $X$ is ordinary and closely related to the multiplication by $p$. It analogously induces a rational map $V_S : S_{X_1} \dashrightarrow S_X$, called the (generalized) Verschiebung, that is generically étale for an ordinary curve and restricts to the morphism $K_{X_1} \to K_X$ induced by $V_J$ on the semi-stable boundary. There are stable vector bundles $E$ over $X_1$ such that $F^*E$ is not semi-stable and the corresponding set  of points in $S_{X_1}$ is precisely the base locus $B(V_S)$ of $V_S$.\\

 In genus 2, $S_X\cong |2\Theta|$ and $K_X$ is a quartic surface, the Kummer surface, with 16 nodes corresponding to the 16 points of $J_X[2]$ and forming a so-called $16_6$ configuration in $|2\Theta|\cong \P^3$. The Verschiebung $V_S$ is given by degree $p$ polynomials. The base locus of Verschiebung is finite and we know from the work of B. Osserman \cite{Oss} (see also \cite{LaP}) that it is reduced for a general curve, of length $ 2(p^3-p)/3$, and that the indeterminacy of $V_S$ can be resolved by a single blowing-up of $S_{X_1}$ at $B(V_S)$. As a conséquence, $V_S$ has degree $(p^3+2p)/3$.

 In genus 2 and characteristic 3, Y. Laszlo and C. Pauly have given a nice description of Verschiebung \cite[Section 6]{LP2}. Consider the family of rank 2 and trivial determinant semi-stable vector bundles parameterized by the blowing-up ${\rm Bl}_2(J_X)$ of $J_X$ at its 2-torsion (viewed as a parameter space for pairs $(L,\,[\varphi])$, $L$ in $J_X$, $[\varphi]$ in $\mathbb{P}H^0(X,\,L^2\tens \omega_X)$) and defined by the exact sequences
 $$ 0 \to \mathcal{E}_{L,\,[\varphi]} \to F_*L \xrightarrow{\varphi} L\tens \omega_{X_1} \to 0$$
 where $L$ also denotes the $p$-twist of $L$ (in such a way that $F^*L\cong L^3$). The moduli property thus induces a morphism ${\rm Bl}_2(J_X) \to S_{X_1}$ that is proved, using fine properties of Kummer surfaces, to factor through an embedding $\alpha : K_X \to S_{X_1}$ such that $V_S$ coincides with the polar map of $\alpha(K_X)$. In particular, there is an equality of divisors
$$V_S^{-1}(K_X)=K_{X_1}+2\alpha(K_{X})$$ in $S_{X_1}$ and  $B(V_S)$ is precisely the singular locus of $\alpha(K_{X})$.

In \cite{Du}, using Prym varieties, we compute  the equations defining $V_S$ in characteristic $p=3,\, 5$ or 7, recovering the results of \cite{LP2}. We can thus check that for a general genus 2 curve, there is a normal and  integral surface $\bar{H}$ of degree $2(p-1)$ in $S_{X_1}$ (the Zariski closure of the set $(V_S^{-1}(K_X))^s$ of stable bundles $E$ such that $F^*E$ is strictly semi-stable) such that the equality of divisors
$$V_S^{-1}(K_X)=K_{X_1}+2\bar{H}$$ in $S_{X_1}$
holds in $S_{X_1}$. Also, we could thus test on many curves (with no exception encountered) the following facts : $B(V_S)$ is contained in the 0-dimensional singular locus ${\rm Sing}(\bar{H})$ of $\bar{H}$, there are singular point of $\bar{H}$ where $V_S$ is defined, these points are mapped onto the singular locus ${\rm Sing}(K_X)$ of $K_X$ and the inverse image of ${\rm Sing}(K_X)$ is 1-dimensional.

In this work, we look at the Zariski closure $\bar{H}$ of the set $(V_S^{-1}(K_X))^s$ of stable rank 2 and trivial determinant bundles $E$ over $X_1$ such that $F^*E$ is non-stable for a general curve of genus $g\geq 2$  in characteristic $p>2(g-1)$. It is dominated by a proper scheme $\widetilde{\mathcal{H}}$ whose study boils down to the study of the Hilbert scheme $\mathcal{H}$ of rank 2 and degree 0 subbundles of $F_*\O_X$. For a point $(E,\,\alpha)$ of $\mathcal{H}$, the quotient fitting in an exact sequence
$$0 \to E \xrightarrow{\alpha} F_*\O_X \to Q_E \to 0$$ is stable and the forgetful map $(E,\,\alpha) \mapsto Q_E$ is a closed immersion if $p>2g$. In this case, $\mathcal{H}$ has dimension at most $2(g-1)$ (its expected dimension as a Hilbert scheme) and  smoothness is therefore given by the vanishing of an $H^1$. The Prym varieties associated to double étale covers of $X$ provide us with (a finite number of) smooth points.

In genus 2, a result of Tong Jilong ensures that the divisor $\Theta_B$ associated to the sheaf of locally exact differential forms is smooth. Pulling-back the canonical exact sequence
\begin{eqnarray}\label{sesdefB1} 0 \to \O_{X_1} \to F_*\O_X \to B \to 0\end{eqnarray}
by an injection $L \to B$ where $L$ is a degree 0 line bundle over $X_1$, provides with essentially every strictly semi-stable points of $\mathcal{H}$.

We derive from these facts our main result :

\begin{theorem}\label{thmbarHintegral} Let $X$ be a general genus 2 curve over an algebraically closed field with characteristic $p\geq 3$. There is an integral and normal surface $\bar{H}$ in $S_{X_1}$, with degree $2(p-1)$, such that the equality $$V_S^{-1}(K_X)=K_{X_1}+2\bar{H}$$ holds scheme-theoretically.\end{theorem}

 All schemes will be (locally) noetherian over a base field $k$. Unless otherwise specified, a point of $k$-scheme $S$ is a closed point (i.e., a $k$-point). We will let ${\bf Sch}/k$ denote the category of schemes over $k$ and ${\bf Sets}$ denote the category of sets.

If $Y$ and $Z$ are schemes over $k$, we let $q_Y$ (resp. $q_Z$) denote the first (resp. the second) projection $Y\times Z \to Y$ (resp. $Y\times Z \to Z$).

When $M$ (resp. $M'$) is a coherent sheaf over a scheme $T$, $\mathcal{E}nd(M)$ (resp. $\mathcal{E}nd_0(M)$, resp. $\mathcal{H}om(M,\,M')$) is the sheaf of germs of $\O_T$-linear endomorphisms of $M$ (resp. traceless endomorphisms of $M$, resp. homomorphisms from $M$ to $M'$). Global sections form the space $\End(M)$ (resp. $\End_0(M)$, resp $\Hom(M,\,M')$). The dual of a vector bundle $M$ over $T$ is denoted by $M^\vee$.

We will say that two elements $f,\,f'$ of $\Hom(M,\,M')$ coincide if they are equal up to multiplication by a non-zero scalar.

\section{Preliminaries}
Let $X$ be a proper and smooth curve of genus $g\geq 2$ over an algebraically closed field $k$ of characteristic $p\neq 2$.

\subsection{Moduli spaces} In the sequel, we will always denote by $S_X$ (resp. $N_X$) the coarse moduli space of rank 2 and trivial determinant (resp. degree 0) semi-stable vector bundles over $X$. These have dimension $4(g-1)-1$ and $5(g-1)$ respectively.

We denote by $S_X^s$ (resp. $N_X^s$) the stable locus of $S_X$ (resp. $N_X$) and by $S_X^{ss}$ (resp. $N_X^{ss}$) the semi-stable boundary of $S_X$ (resp. $N_X$).  The semi-stable boundary of $S_X$ (resp. of $N_X$) identifies with the Kummer variety $K_X$ of $X$ and it is precisely the singular locus of $S_X$ (resp. of $N_X$) unless $g=2$. Whenever a scheme $T$ parameterizes a family $\mathcal{E}$ of rank 2 and trivial determinant (resp. degree 0) over $X$, we let $T^s$ denote the inverse image of $S_X^s$ (resp. $N_X^s$) by the map $T \to S_X$ (resp. $T \to N_X$) induced by the moduli property. Similarly, $T^{ss}$ will denote the inverse image of $S_X^{ss}$ (resp. $N_X^{ss}$).

It is well-known that $S_X$ and $N_X$ are normal (even smooth when $g=2$) and projective varieties over $k$ and that there is a finite étale morphism $t : S_{X}\times
J_{X} \to N_{X}$ (mapping $(E,\,L)$ to $E\tens L$)
of degree $2^{2g}$. Taking the determinant of a vector bundle gives a morphism $\det : N_X \to J_X$ and letting $[2] : J_X \to J_X$ denote multiplication by 2, the diagram
\begin{equation}\label{diagNXSXJX}
\unitlength=0.6cm
\begin{picture}(8,4)
\put(0.8,0){$N_X$} \put(0,3){$S_X\times J_X$} \put(6,0){$J_X$} \put(6,3){$J_X$}
\put(3,3.2){\vector(1,0){2.6}} \put(2,0.2){\vector(1,0){3.6}}
\put(1.4,2.8){\vector(0,-1){2}} \put(6.4,2.8){\vector(0,-1){2}}
{\small \put(0.8,1.5){$t$} \put(6.8,1.5){$[2]$} \put(4,3.5){$q_J$} \put(3.4,0.5){$\det$}}
\end{picture}
\end{equation}
is cartesian. Recall that, at a stable point $E$ of $N_{X}$ (resp. $S_X$), the tangent space $T_EN_X$ (resp. $T_ES_X$) at $E$ is canonically isomorphic to $H^1(\mathcal{E}nd (E))$ (resp. $H^1(\mathcal{E}nd_0 (E))$) and the tangent map of $\det : N_X \to J_X$ is the trace map ${\rm Tr} : H^1(\mathcal{E}nd (E)) \to H^1(\O_X)$.

When $k$ has characteristic $p$, we let $M_{X}^Q$ denote the moduli space $$M_{X}(p-2,\,(p-1)(g-1))$$ of rank $p-2$ and degree $(p-1)(g-1)$ semi-stable vector bundles over $X$. It is proper and normal of dimension $(p-2)^2(g-1)-1$ and its tangent space at a stable point $Q$ is canonically isomorphic to $H^1(\mathcal{E}nd (Q))$.
\subsection{An upper bound for slopes of rank $r$ subbundles of $F_*\eta$}

We apply a construction of \cite{JRXY} relying on adjunction and a theorem of P. Cartier (see \cite[Theorem 5.1]{Ka}) to prove the following key-lemma

\begin{lemma} \label{lemsubF*xi} Let $X$ be a genus $g\geq 2$ curve over $k$ and  let $\eta$ be a degree $d$ line bundle over $X$.
For all $1\leq r\leq p$, a rank $r$
subbundle $M$ of $F_*\eta$ has slope
$$\mu\leq \mu_r(\eta)=\frac{(r-1)(g-1)+d}{p}$$
Furthermore, there is equality if and only if $F^*M$ has a filtration
$$\{0\}=M_r\subset M_{r-1}\subset ... \subset M_1 \subset M_0 =F^*M$$ with $M_i/M_{i+1}\cong \omega_X^{\tens i}\tens \eta$ for all $0 \leq i\leq r-1$.
\end{lemma}
\begin{proof} The case $r=1$ comes directly from adjunction formula. If $M$ is a  rank $r\geq 2$ subbundle of $F_*L$, we write $M_0=F^*M$ and we let $M_1$ denote the kernel of the non zero map  $\varphi_0 : M_0=F^*M \to \eta$ corresponding to $M \hookrightarrow F_*\eta$ by adjunction. There is an exact sequence
$$0 \to M_1 \to M_0 \to M_0/M_1 \to 0$$
over $X$ and $\deg (M_0/M_1)\leq d$  with equality if and only if $M_0/M_1\cong \eta$.
The second
fundamental form of the Cartier's connection $\nabla : M_0 \to M_0\tens \omega_X$ associated to this exact sequence is a $\O_X$-linear map $\varphi_1 : M_1 \to
(M_0/M_1) \tens \omega_X$ that is not zero for otherwise $M_1$ descends to a subbundle of $M_0$, hence of $F_*\eta$, and the composite $M_1 \to M_0 \xrightarrow{\varphi_0} \O_X$ does not vanish by functoriality of adjunction.

If $r=2$, $M_1$ is a line bundle and one has $\deg M_1 \leq \deg (M_0/M_1) +2(g-1) \leq 2(g-1)+d$ with equality if and only if $M_1\cong M_0/M_1\tens \omega_X$, hence $$\mu=\frac{1}{p}\mu(M_0)=\frac{1}{p}\frac{\deg (M_1)+\deg(M_0/M_1)}{p}\leq \mu_2(\eta)$$ with equality if and only if $M_0/M_1\cong \eta$ and $M_1\cong \omega_X \tens \eta$.

If $r\geq 3$, we let $M_2$ denote the kernel of $\varphi_1$, a rank $r-2$ subbundle of $M_0$ and we let $\varphi_2 : M_2 \to M_0/M_2\tens \omega_X$ denote the second fundamental form of $\nabla$ associated to the exact sequence
$$0 \to M_2  \to M_0 \to M_0/M_2 \to 0$$
For the same reason as above, $\varphi_2$ does not vanish and we construct inductively a filtration
$$\{0\} =M_r \subset M_{r-1} \subset ... \subset M_1 \subset M_0=F^*M$$
of $M_0$ and a family of $\O_X$-linear maps $$\varphi_k : M_k \to (M_0/M_k)\tens \omega_{X}$$
where $$M_{k}=\ker (\varphi_{k-1} : M_{k-1} \to (M_0/M_{k-1})\tens \omega_X)$$
and $\varphi_k$ is the second fundamental of $\nabla$ associated to
the exact sequence
$$0 \to M_k \to M_0 \to M_0/M_k \to 0$$
By construction, the map $\varphi_k$ induces an injection $(M_{k}/M_{k+1}) \hookrightarrow (M_{k-1}/M_{k})\tens \omega_X$ which yields by induction an inequality
$$\mu=\frac{1}{p}\mu(M_0)=\frac{1}{pr}\sum_{k=0}^{r-1}\deg(M_k/M_{k+1})\leq \frac{1}{pr}\sum_{k=0}^{r-1}(2(g-1)k+d)=\mu_r(\eta)$$ with equality if and only if $M_{k}/M_{k+1}\cong \omega_X^{k}\tens \eta$ for all $0\leq k\leq r-1$.
 \end{proof}

\begin{remark} \label{rkquotF*xi} It immediately follows that a rank $r$ quotient sheaf $M'$ of $F_*\eta$ has slope
$$\mu\geq \nu_r(\eta)=\frac{(2p-r-1)(g-1)+d}{p}=\mu(F_*\eta)+\frac{(p-r)(g-1)}{p}$$
\end{remark}

In the sequel, we simply write $\mu_r$ (resp. $\nu_r$) for $\mu_r(\O_X)$ (res. $\nu_r(\O_X)$).

\begin{remark}It follows from Lemma \ref{lemsubF*xi} that $F_*\eta$ is stable. When one tries to provide an upper bound $\mu_r(W)$ for rank $r$ subundles of $F_*W$ where $W$ is a stable rank $n$ bundle over $X$, one can proceed analogously to define a filtration $$\{0\} =M_r \subset M_{r-1} \subset ... \subset M_1 \subset M_0=F^*M$$ in the same way as above. In particular, there is a series of inclusions
$$M_{r-1} \subset (M_{r-2}/M_{r-1})\tens \omega_X \subset ... \subset (M_1/M_2 ) \tens \omega_{X}^{r-2} \subset (M_0/M_1) \tens \omega_{X}^{r-1} \subset W\tens \omega_{X}^{r-1}$$
However, since we are dealing with higher rank vector bundles, one is confronted to the fact that a non zero map needs not be of maximal rank. Still, if ${\rm rk}(M_k/M_{k+1})<{\rm rk}(M_{k-1}/M_{k})$, there is a quotient bundle $W\to W''$ with kernel $W'$ such that the composite $$M_k/M_{k+1} \subset (M_0/M_1) \tens \omega_{X}^{k} \subset W\tens \omega_{X}^{k} \twoheadrightarrow W''\tens \omega_X^k$$ is zero. Letting $M'$ denote the kernel of the composite $M \subset F_*W \twoheadrightarrow F_*W''$, this is a subbundle of $F_*W'$ and we can analogously define a sequence
 $$\{0\} =M'_{r'} \subset M'_{r'-1} \subset ... \subset M'_1 \subset M'_0=F^*M'$$
The point is then to notice that $M'_i=M_i$ for all $i\geq k$ and proceed by induction on $n$ to prove that $r\leq p$. The stability of $W$ thus induces the stability of $F_*W$ (see \cite{Sun} for an another proof much in the same spirit and a  generalization of the argument for higher dimensional varieties).\end{remark}

\subsection{Recollections on the divisor $\Theta_B$}\label{subsecprelimB}
Assume that $k$ has characteristic $p\geq 3$ and recall that the differential $d : \O_X \to \omega_X$ kills all $p$-powers and induces the two short exact sequences \eqref{sesdefB1} and
\begin{equation} \label{sesdefB2} 0 \to B \to F_*\omega_X \xrightarrow{c} \omega_{X_1} \to 0\end{equation} of vector bundles over $X_1$, where $c$ is the Cartier's operator (see \cite[Section 4]{Ra}). The bundle $B$ of locally exact differential forms has rank $p-1$ and slope $g-1$. Taking $\eta=\O_X$ in the lemma \ref{lemsubF*xi}, we find that $B$ is stable (see\cite{Jo} for another argument) and it provides a filtration \begin{equation}\label{eqfiltF*F*OX}
\{0\}=B_p\subset B_{p-1}\subset .... \subset B_1=F^*B\subset B_0= F^*F_*\O_X
\end{equation}
with $B_i/B_{i+1}\cong \omega_X^{\tens i}$ for all $0\leq i\leq p-1$.
(recall that the filtration \eqref{eqfiltF*F*OX} of $B_1$ was obtained in \cite[Section 4]{Ra} using the interpretation of $B_1$ as the augmentation ideal of the $(p-1)$-th thickening of the diagonal in $X\times _kX$).

Recall
that the bundle $B$ is endowed with a non-degenerate skew-symmetric bilinear form
\begin{equation}\label{eqBtensBomega} \Upsilon : B \tens B \to \omega_{X_1} \ \text{ (or, equivalently, }\bar{\Upsilon} : \Lambda^2B \to \omega_{X_1})\end{equation}
 Because $F^*F_*\O_X=B_0\cong \O_X\oplus B_1$, there is an isomorphism
$$F_*\O_X \tens F_*\O_X \cong F_*\O_X \oplus F_
*\O_X \tens B$$ where the projection on the first factor is the canonical multiplication map. Thus, the subbundle $\Lambda^2F_*\O_X \subset F_*\O_X \tens F_*\O_X$ actually lies in  $F_*\O_X \tens B$. Also, it is tautological that the composite $\Lambda^2F_*\O_X \to F_*\O_X \tens B \to B \tens B$ factors through $\Lambda^2B \to B \tens B$ and since there is an exact sequence
$$0 \to  B \to \Lambda^2F_*\O_X \to \Lambda^2B \to 0$$
we find that the extension class of $\Lambda^2F_*\O_X$ in $\Ext^1(\Lambda^2B,\,B)$ is the image of the extension class of $F_*\O_X $ in $\Ext^1(B,\,\O_X)$ via the map
$$\Ext^1(B,\,\O_X) \xrightarrow {1_B\tens 1} \Ext^1(B\tens B,\,B)\xrightarrow {\Lambda^2B \to B \tens B} \Ext^1(\Lambda^2B,\,B)$$ Now, consider the exact sequence \eqref{sesdefB2} and pull-it back by $\bar{\Upsilon} : \Lambda^2B\to \omega_{X_1}$. Because the extension class of $F_*\omega_X$ in $\Ext^1(\omega_{X_1},\,B)$ is that of $F_*\O_X$ in $\Ext^1(B,\,\O_X)$ via the isomorphism $B \xrightarrow \sim B^\vee \tens \omega_{X_1}$ deduced from $\Upsilon$, there is a map
\begin{equation}\label{eqdefwDelta} \widetilde{\Upsilon} : \Lambda^2F_*\O_X \to F_*\omega_X \end{equation}
fitting in the commutative diagram with exact rows and columns
\begin{equation}\label{diagLambda2FODelta}\begin{array}{ccccccc}
 & & & 0& & 0\\
 & & & \downarrow & & \downarrow &\\
  & & & \Gamma & = & \Gamma & \\
   & & & \downarrow & & \downarrow &\\
 0 \to & B & \to & \Lambda^2F_*\O_X & \to & \Lambda^2 B & \to 0\\
 & || & & \downarrow & & \downarrow \\
 0 \to & B & \to & F_*\omega_X & \to & \omega_{X_1} & \to 0\\
 & & & \downarrow & & \downarrow &\\
 & & & 0& & 0\end{array}\end{equation}
where $\Gamma$ is the kernel of $\bar{\Upsilon}$.

In \cite{Ra}, M. Raynaud proves that the set $$\{L\in J_{X_1}|\Hom_{X_1}(L,\,B)\neq 0\}$$ is the support of a well-defined divisor $\Theta_B$. The latter is algebraically
equivalent to $(p-1)\Theta$ and totally symmetric (in the sense of Mumford), hence comes from a divisor on the Kummer variety $K_{X_1}$ of $X_1$. In case $X$ is an ordinary curve, $\Theta_B$ does not go through the origin of $J_{X_1}$ (in particular $H^0(B)=H^1(B)=0$). We recall a few results on $\Theta_B$ that Jilong Tong, a former student of M. Raynaud, obtains by degeneration technics (PhD Thesis \cite{To}). Let us first set the following

\begin{definition}\label{defHilb10B} Let $\mathcal{H}_B$ be the {\em Quot scheme
representing the functor} $\mathcal{H}ilb_{1,\,0}(B)$ from the
category of schemes over $k$ to {\bf Sets} and defined by $$T \mapsto
\mathcal{H}ilb_{1,\,0}(B)(T)=\left\{\begin{array}{c}
\text{ rank } p-2 \text{ quotient sheaves } q_{X_1}^*B \twoheadrightarrow \mathcal{Q} \\
\text{ over
}X_1\times T, \text{ flat over
} T, \\
\text{ with } \deg
\mathcal{Q}_t=p-1 \ \forall t \in T(k)
\end{array}\right\}/\sim $$\\
The universal line subbundle over
$X_1\times \mathcal{H}_B$ induces a morphism $\mathcal{H}_B \to \Theta_B \subset J_{X_1}$ that can be seen as a "projective" bundle over $\Theta_B$ (meaning that the fiber above a point $L$ of $\Theta_B$ is isomorphic to $\mathbb{P}\,\Hom_{X_1}(L,\,B)^\vee$).
\end{definition}

\begin{theorem}[\cite{To}] \label{thmTongJ} Let $X$ be a general proper and smooth curve of genus $g\geq 2$ over an algebraically closed field $k$.\\
(1) If $g=2$, the divisor $\Theta_B$ is smooth and the map $\mathcal{H}_B \to \Theta_B$ is an isomorphism (Theorem 3.2.2, p. 60).\\
(2) If $g\geq 3$, the divisor $\Theta_B$ is normal (Theorem 3.2.3, p. 60).\\
(3) If $g\geq 2$, $J_{X_1}[p\,]\setminus \{0\}$ is the whole set of points of finite order in $\Theta_B$ and $\Theta_B$ is smooth at each of them (Proposition 3.3.2.1, p. 63). In particular, if $L$ is an order $p$ line bundle over $X_1$, $\dim \Hom(L,\,B))=1$ (Corollary 2.3.2, p.50).
\end{theorem}

\section{Rank 2 and degree 0 subbundles of $F_*\O_X$ for $g\geq 2$} From now on, we assume that the base field $k$ has characteristic $p>2(g-1)$ and that
the curve $X$ is general. In particular, it is ordinary and we have a finite étale group scheme \begin{equation} \label{eqdefG} G:=J_{X_1}[p\,]_{\rm red} \cong \ker(V_J : J_{X_1} \to J_X)\cong (\mathbb{Z}/p\,\mathbb{Z})^g\end{equation}
(it is well known that this statement holds for any ordinary abelian variety of dimension $g\geq 1$).

\subsection{Some definitions}\label{subsecsomedef}

Our main object of study will be the following
\begin{definition}\label{defHeta} Let $\eta$ be a degree $d$ line bundle over $X$ with $d\leq 0$ and let $\mathcal{H}_\eta$ be the {\em Quot-scheme
representing the functor} $\mathcal{H}ilb_{2,\,0}(F_*
\eta)$ from ${\bf Sch}/k$ to {\bf Sets} and defined by $$T \mapsto
\mathcal{H}ilb_{2,\,0}(F_*\eta)(T)=\left\{\begin{array}{c}
\text{rank }p-2\text{ quotient sheaves } q_{X_1}^*(F_*\eta) \twoheadrightarrow \mathcal{Q} \\ \text{ over
}X_1\times T, \text{ flat over }T, \\
\text{ with }\deg
\mathcal{Q}_t=(g-1)(p-1)+d \ \forall t \in T(k)
\end{array}\right\}/\sim $$
where $q_{X_1} : X_1\times T \to X_1$ denotes the first projection.\\
If $\eta=\O_X$, we let $\mathcal{H}$ denote the Quot-scheme $\mathcal{H}_{\O_X}$.
\end{definition}

For a closed point $(E,\,\alpha)$ of $\mathcal{H}$, we let $Q_E$ denote the quotient sheaf fitting in the exact sequence
\begin{equation}\label{sesdefQ}0 \to E \xrightarrow{\alpha} F_*\O_X \to Q_E \to 0.\end{equation}
By adjunction, there is a non-zero map ${\rm ad}(\alpha) : F^*E \to \O_X$ that factors through an injection $s_\Delta : \O(-\Delta) \hookrightarrow \O_X$, where $\Delta$ is an effective (and possibly trivial) divisor. If $\Delta$ is has maximal degree among those with this property, the map $F^*E \to \O(-\Delta)$ is surjective and there is an exact sequence
\begin{equation}\label{sesF*EDelta}0 \to \O(\Delta)\tens F^*(\det E) \to F^*E \to \O(-\Delta) \to 0\end{equation}One checks that the linear equivalence class $\O(\Delta)$ of $\Delta$ is uniquely defined and in case $\O(\Delta)$ is non-trivial, we call its degree the \emph{degree of Frobenius-destabilization} of $E$.

In application to the Hilbert scheme $\mathcal{H}$, the lemma \ref{lemsubF*xi} has the following immediate consequence

\begin{corollary} \label{cordegDelta<g} Assume that $(E,\,\alpha)$ is a point of $\mathcal{H}$ and that $F^*E$ fits in an exact sequence as above. The second fundamental form
\begin{equation*}\varphi_1 : \O(\Delta)\tens F^*(\det E) \to \O(-\Delta) \tens \omega_X\end{equation*} of Cartier's  connection  associated to this exact sequence does not vanish.\\ In particular, $H^0(X,\,\O(K-2\Delta)\tens \det (F^*E))\neq 0$ and $0\leq \deg(\Delta)\leq g-1$.\end{corollary}

If $\Delta=0$, $F^*E$ is semi-stable and one has $\Hom_{X_1}(E,\,F_*\O_X) \cong \Hom_X(F^*E,\,\O_X)$ by adjunction. We define  $$\mathcal{H}_2
 :=\{(E,\,\alpha) \text{ in }\mathcal{H}|\, F^*E\cong \O_X\oplus \O_X\}.$$

\begin{remark}\label{rkecupphi1=0} Assume that $F^*E$ is semi-stable and let $e$ in $H^1(X,\,F^*(\det E))$ denote the extension class of \eqref{sesF*EDelta}. Letting $\varphi_1$ denote as above the second fundamental form of Cartier's connection associated to this exact sequence, a straightforward cocycle computation yields that the cup-product $e\, \cup \varphi_1$ in $H^1(\omega_X)$ vanishes.
\end{remark}

If $\Delta$ is non trivial, it follows from \cite[Satz 2.5]{LaS} that the extension \eqref{sesF*EDelta} is non-split and one has  $\Hom_{X_1}(E,\,F_*\O_X) \cong H^0(X,\,\O(\Delta))$.
The universal inclusion $s : \O(-\Delta) \to \O_X$ over $X\times |\Delta|$ gives rise to an inclusion $F_*\O(-\Delta) \xrightarrow {F_*s} F_*\O_X$ over $X_1\times |\Delta|$ and the universal property of $\mathcal{H}$ yields a morphism
$$\mathcal{H}_{\O(-\Delta)} \times |\Delta| \to \mathcal{H}$$One checks that this is a closed immersion and we denote $\mathcal{H}_{|\Delta|}$ its image in $\mathcal{H}$.
We check similarly that if $\Delta$ and $\Delta'$ are two effective divisors with $\Delta-\Delta'$ effective, there is a series of inclusions $\mathcal{H}_{|\Delta|}\subset \mathcal{H}_{|\Delta'|} \subset \mathcal{H}$.\\

Recall the generalized Nagata-Segre theorem (relying on a result of Hirchowitz \cite{Hir}, see also \cite[Theorem 2.1]{LaP}) :
\begin{theorem} For any rank $r$ and degree $\delta$ vector bundle $W$ over any smooth curve $X$ of genus $g$ and for any integer $n$ with $1\leq n\leq r-1$, there is a subbunble $E\subset W$ such that the inequality
$$\mu(E)\geq \mu(W)-\frac{r-n}{r}(g-1)-\frac{\epsilon}{rn}$$ holds, where $\epsilon$ is the only integer with $0\leq \epsilon \leq r-1$ and $\epsilon +n(r-n)(g-1)\equiv n\delta$ mod $r$.\end{theorem}

\begin{corollary}\label{corHetanonempty} For any line bundle $\eta$ with degree $1-g\leq d\leq 0$, the Quot-scheme $\mathcal{H}_\eta$ is non empty.\end{corollary}
\begin{proof}Taking $W=F_*\eta$ and $n=2$, we obtain $\epsilon\equiv 2(g-1+d)$ and considering our hypothesis on $p,\ g$ and $d$, $\epsilon=2(g-1+d)$. Now, the theorem provides us with a rank 2 subbundle of $F_*\eta$ with slope 0.\end{proof}

\bigskip

On the one hand, the universal subbundle $\mathcal{E}_\eta$ in
the universal exact sequence
\begin{equation*}
0 \to \mathcal{E}_\eta \xrightarrow \alpha q_{X_1}^* (F_*\eta) \xrightarrow \beta \mathcal{Q}_\eta \to 0\end{equation*} over
$X_1\times \mathcal{H}_\eta$ obviously induces a forgetful morphism \begin{equation}\label{defh_E} h_{\mathcal{E}_\eta} : \mathcal{H}_\eta \to N_{X_1}\end{equation} via the moduli property.
Using the morphism $t : S_{X_1}\times J_{X_1} \to N_{X_1}$, we can consider the fiber product $\mathcal{H}_\eta \times _{N_{X_1}} (S_{X_1}\times J_{X_1})$. The map $\mathcal{H}_\eta \times _{N_{X_1}} (S_{X_1}\times J_{X_1})\to \mathcal{H}_\eta$ is étale of degree $2^{2g}$ and there is a map $\mathcal{H}_\eta \times _{N_{X_1}} (S_{X_1}\times J_{X_1}) \to S_{X_1}\times J_{X_1}$ above $h_{\mathcal{E}_\eta}$.

On the other hand, the group $G$ naturally acts upon $\mathcal{H}_\eta$ for any $\eta$. Indeed, for any $\gamma$ in $G$, use the canonical isomorphism $F_*\eta \tens \gamma \xrightarrow \sim  F_*\eta$ and define $\tau_\gamma$ as the automorphism of $\mathcal{H}_\eta$ defined by \begin{equation}\label{deftauN}\tau_\gamma(\mathcal{E}_T) \cong \mathcal{E}_T\tens q_{X_1}^*(\gamma)\end{equation} for any such $\gamma$ and for any $T$-point of $\mathcal{H}_\eta$. This action is obviously free and it defines an étale quotient map $\mathcal{H}_\eta \to \mathcal{H}_\eta/G$ of degree $p^g$. This action can be lifted to a free action of $G$ upon $\mathcal{H}_\eta \times _{N_{X_1}} (S_{X_1}\times J_{X_1}) $ and the quotient map is similarly étale finite of degree $p^g$.

\begin{definition} \label{defwHeta} We define
\begin{equation*}\wH_\eta:=\left[\mathcal{H}_\eta \times _{N_{X_1}} (S_{X_1}\times J_{X_1})\right]/G \ \text{ (resp. }\wH:=\left[\mathcal{H} \times _{N_{X_1}} (S_{X_1}\times J_{X_1})\right]/G).\end{equation*}
We let $\wH_2$ denote the set of points in $\wH$ such that $F^*E$ is split isotypic.\\
We let $\wH_{|\Delta|}$ denote the closed subscheme of $\wH$ corresponding to $\mathcal{H}_{|\Delta|}$ for any effective divisors $\Delta$ as above. It is the image of the induced closed immersion $\wH_{\O(-\Delta)}\times |\Delta| \hookrightarrow \wH$.

A point in $\wH_{|\Delta|}$ thus can be seen as a triple $(E,\,\xi,\,\alpha)$ where $E$ is a semi-stable rank 2 vector bundle with trivial determinant, where $\xi$ is a degree 0 line bundle over $X$ and where $\alpha$, defined up to a multiplicative scalar, is an embedding $E\hookrightarrow F_*(\O(-\Delta)\tens \xi)$.
\end{definition}
\begin{definition}\label{defwh+barH}There is an induced map
$\wH_\eta \to S_{X_1}\times J_{X}$ and we define \begin{equation*}\widetilde{h}_\eta :  \wH_\eta \to S_{X_1} \text{ and } \widetilde{V}_\eta : \wH_\eta \to J_X \text{ (resp. }\widetilde{h}:  \wH \to S_{X_1} \text{ and } \widetilde{V} : \wH \to J_X) \end{equation*} as its composites with the projections.\\
We let \begin{equation*} \bar{H}:={\rm im} (\widetilde{h} :  \wH \to S_{X_1})\end{equation*} denote the \emph{scheme-theoretic} image of $\wH$ ($\bar{H}^s$ and $\bar{H}^{ss}$ will denote respectively its stable and its strictly semi-stable loci) and
\begin{equation*}\bar{H}_2:=\{E \text{ in } \bar{H}^s| F^*E\cong \tau\oplus \tau,\ \tau \text{ in }J_X[2]\}\end{equation*}\end{definition}

\begin{lemma} \label{lemfib(wh)} Assume that $E$ is a stable vector bundle over $X_1$ with rank 2 and trivial determinant (viewed as a point of $S_{X_1}$). Then,\\
(1) the bundle $E$ lies in $\bar{H}$ if and only if $F^*E$ is not stable .\\
(2) the bundle $E$ lies in the image of $\wH_{|\Delta|}$ for some non trivial $\Delta$ if and only if $F^*E$ is unstable.\\
(3) If $F^*E$ is semi-stable, then $\widetilde{h}^{-1}(E)$ contains at most 2 (opposite) points unless $E$ lies in $\bar{H}_2$. In this case, $\widetilde{h}^{-1}(E)\cong \P^1$ with $\widetilde{V}(\widetilde{h}^{-1}(E)) \subset J_X[2]$.\\
 In any case, the image of $\widetilde{V}(E,\,\xi,\,\alpha)$ through the composite $J_X \to K_X \subset S_X$ is $V_S(E)$.\\
(4) If $F^*E$ is unstable, the fibre $\widetilde{h}^{-1}(E)$ identifies with the set ${\rm Sym}^d(X)$ of degree $d$ effective divisors over $X$, where $d$ is degree of Frobenius-destabilization of $E$.\\
 The map $\widetilde{V}$ restricts on the fibre $\widetilde{h}^{-1}(E)$  to the usual map ${\rm Sym}^d(X) \to J_X$ (defined up to translation).
\end{lemma}
\begin{proof}(1), (2) and (3) immediately follow from the adjunction formula.

As for (4), assume that $E$ is a stable point of $\bar{H}$ lying in the image of $\wH_{|\Delta|}$. A point $(E,\,\xi,\,\alpha)$ lies in the fibre $h^{-1}(E)$ if and only if $F^*E$ fits in an exact sequence
$$0 \to \O(\Delta)\tens \xi^{-1} \tens \to F^*E \to \xi \tens \O(-\Delta) \to 0$$ for some non-trivial effective divisor $\Delta$. Then, for any effective divisor $\Delta'$ of degree $d$, there is a surjective map $$F^*E\tens \O(\Delta-\Delta') \to \xi\tens  \O(-\Delta')$$ giving rise to an injective map $E \hookrightarrow F^*(\xi \tens \O(-\Delta')) \xrightarrow{F_*(1\tens s_{\Delta'})} F_*\xi$. Calling this map $\alpha'$, $(E,\,\xi\tens \O(\Delta-\Delta'),\,\alpha')$ is a point of $h^{-1}(E)$. Notice that if $\Delta$ and $\Delta'$ are two distinct linearly equivalent divisors, the map $\alpha$ and $\alpha'$ do not coincide.

Conversely, if $(E,\,\xi',\,\alpha')$ is an other point in $h^{-1}(E)$
fitting in an exact sequence $$0 \to \O(\Delta')\tens \xi'^{-1} \tens \to F^*E \to \xi' \tens \O(-\Delta') \to 0$$
the composite $\O(\Delta)\tens \xi^{-1} \tens \to F^*E \to \xi' \tens \O(-\Delta')$ is necessarily 0 and $\O(\Delta)\tens \xi^{-1}$ is a line subbundle of $\O(\Delta')\tens \xi'^{-1}$. By symmetry, one obtains $\O(\Delta)\tens \xi^{-1}\cong \O(\Delta')\tens \xi'^{-1}$. Since the map $\widetilde{V}$ is just the forgetful map $(E,\,\xi,\,\alpha) \mapsto \xi$, this concludes the proof.\end{proof}

\subsection{The "second fundamental form" map}

In the sequel, we use the notations of Subsection \ref{subsecprelimB}. Recall that we let $B_0$ (resp. $B_1$) denote $F^*F_*\O_X$ (resp. $F^*B$) and $\psi_0 : B_0 \to \O_X$ (resp. $\psi_1 : B_1 \to \omega_X$) denote the evaluation map  (resp. the second fundamental form of Cartier's connection $\nabla : B_0 \to B_0 \tens \omega_X$ associated to the exact sequence $0 \to B_1 \to B_0 \xrightarrow{\psi_0} \O_X \to 0$). Notice that the adjoint of $\psi_1$ is the injection $B \hookrightarrow F_*\omega_X$ of the exact sequence \eqref{sesdefB2}. Eventually, we will use the map $\widetilde{\Upsilon} : \Lambda^2F_*\O_X \to F_*\omega_X$. Notice that pulling back by Frobenius splits the middle row in the diagram \eqref{diagLambda2FODelta} and that the map ${\rm ad}(\widetilde{\Upsilon}) : F^*(\Lambda^2F^*\O_X) \to \omega_X$ is nothing but the projection on $B_1$ followed by $\psi_1 : B_1 \to \omega_X$.

\begin{proposition} \label{lemPsi_1=adLambda2}Let $(E,\,\alpha)$ be a point of $\mathcal{H}$ and assume $F^*E$ fits in an exact sequence \eqref{sesF*EDelta}. The adjoint ${\rm ad}(\Phi_1(E)) : F^*(\det E) \to \omega_X$ of the composite
$$\Phi_1(E) : \Lambda^2E\xrightarrow{\Lambda^2\alpha} \Lambda^2F_*\O_X \xrightarrow{\widetilde{\Upsilon}} F_*\omega_X$$ coincides (up to a multiplicative scalar) with the composite
$$F^*(\det E) \xrightarrow{s_\Delta} F^*(\det E) \tens \O(\Delta) \xrightarrow {\varphi_1} \O(-\Delta) \xrightarrow{s_\Delta} \omega_X$$
\end{proposition}
\begin{proof}Take a point $(E,\,\alpha)$ in $\mathcal{H}$ and pull back the exact sequence \eqref{sesdefQ} back to $X$ by Frobenius. Using the notations of the proof of Lemma \ref{lemsubF*xi}, let $E_0$ (resp. $\varphi_0$) denote the Frobenius inverse image of $E$ (resp. the adjoint of $\alpha$). By the very definition of adjunction, $\varphi_0$ is the composite $E_0 \xrightarrow {F^*\alpha} F^*F_*\O_X \xrightarrow {\psi_0}\O_X$.

Since $F^*\alpha$ induces an injective map
 $\rho_1 : \O(\Delta)\tens F^*(\det E) \to B_1$ and since $F^*\alpha$ commutes with the canonical connections on $E_0$ and $B_0$, the composites
\begin{equation}\label{eqrho1psi1}\O(\Delta)\tens \det E_0 \xrightarrow {\rho_1}B_1 \xrightarrow{\psi_1} \omega_X \text{ and }\O(\Delta)\tens \det E_0 \xrightarrow {\varphi_1} \O(-\Delta)\tens \omega_X \xrightarrow{s_\Delta}  \omega_X\end{equation} coincide. Taking Frobenius inverse image and second exterior powers commute and $F^*(\Lambda^2\alpha)$ agrees with the second exterior power of $$F^*\alpha =(\varphi_0,\,F^*\alpha_B) : E_0 \to \O_X\oplus B_1$$ where $\alpha_B$ denotes the composite map $E\xrightarrow \alpha F_*\O_X \to B$. A local computation yields the decomposition
$$\begin{array}{c}F^*(\Lambda^2\alpha)=(\widetilde{\rho}_1,\,F^*(\Lambda^2\alpha_B)) : \det E_0 \to B_1\oplus \Lambda^2B_1 \end{array}$$
where $\widetilde{\rho}_1$ is the composite $\det E_0 \xrightarrow{s_\Delta} \O(\Delta)\tens \det E_0 \xrightarrow{\rho_1} B_1$.
Noticing that the adjoint of $\Phi_1(E)$ is the composite
$$\det E_0=F^*(\Lambda^2E) \xrightarrow{F^*(\Lambda^2\alpha)}F^*(\Lambda^2F^*\O_X) \xrightarrow{{\rm ad}(\widetilde{\Upsilon})} \omega_X$$ concludes the proof.\end{proof}

\bigskip

If one writes a line bundle $L$ over $X$ under the form $\omega_X\tens \O(-x_1-...-x_{2g-2})$, where the $x_i$'s are points of $X$, one has
$$\Hom(L,\omega_X)\cong H^0(\O(x_1+...+x_{2g-2}))$$
Consider the morphism $${\rm Sym}^{2(g-1)}(X) \to J_X$$ defined by
$x_1+...+x_{2g-2} \mapsto \omega_X\tens \O(-x_1-...-x_{2g-2})$ and call $\mathcal{B}(J_{X_1})$ and $\mathcal{B}_2(J_{X})$ respectively the cartesian products
\begin{equation*}
\unitlength=0.6cm
\begin{picture}(20,4)
\put(0.8,0){$J_{X_1}$} \put(0.5,3){$\mathcal{B}(J_{X_1})$} \put(6,0){$J_X$} \put(4.8,3){${\rm Sym}^{2(g-1)}(X)$}
\put(2.6,3.2){\vector(1,0){1.8}} \put(2,0.2){\vector(1,0){3.6}}
\put(1.4,2.8){\vector(0,-1){2}} \put(6.4,2.8){\vector(0,-1){2}}
{\small \put(3.4,0.5){$V_J$}}
\put(9,1.4){\text{ and }}
\put(12,0){\put(0.8,0){$J_{X}$} \put(0.5,3){$\mathcal{B}_2(J_{X})$} \put(6,0){$J_X$} \put(4.8,3){${\rm Sym}^{2(g-1)}(X)$}
\put(2.6,3.2){\vector(1,0){1.8}} \put(2,0.2){\vector(1,0){3.6}}
\put(1.4,2.8){\vector(0,-1){2}} \put(6.4,2.8){\vector(0,-1){2}}
{\small \put(3.4,0.5){$[2]$}}}
\end{picture}
\end{equation*}
It is clear that $\mathcal{B}(J_{X_1})$ represents the functor from ${\bf Sch}/k$ to {\bf Sets} defined by $$T \mapsto
\mathcal{H}ilb_{1,\,0}(F_*\omega_X)(T)=\left\{\begin{array}{c}
\text{rank }p-1\text{ quotient sheaves } q_{X_1}^*(F_*\omega_X) \twoheadrightarrow \mathcal{R} \\ \text{ over
}X_1\times T, \text{ flat over }T, \\
\text{ with }\deg
\mathcal{R}_t=(p+1)(g-1)\ \forall t \in T(k)
\end{array}\right\}/\sim $$
and since any $T$-point $(\mathcal{E}_T,\alpha_T)$ of $\mathcal{H}$ induces a morphism
$$\det \mathcal{E}_T \xrightarrow{\Lambda^2\alpha_T} q_{X_1}^*(\Lambda^2F_*\O_X) \xrightarrow{q_{X_1}^*\widetilde{\Delta}} q_{X_1}^*(F_*\omega_X)$$
that does not vanish at any closed points of $T$, there are  natural morphisms
$$\Phi_1 : \mathcal{H} \to \mathcal{B}(J_{X_1}) \ \text{ and }\ \widetilde{\Phi}_1 : \wH \to \mathcal{B}_2(J_{X})$$
We will refer to these maps
as the "second fundamental form" maps.

It is classical that ${\rm Sym}^{2(g-1)}(X)$ is smooth over $k$, of dimension $2(g-1)$ and so are $\mathcal{B}(J_{X_1})$ and $\mathcal{B}_2(J_{X})$. In particular, in genus 2, the map ${\rm Sym}^2(X) \to J_X$ (resp. $\mathcal{B}_2(J_X)$) agrees with the blowing-up of $J_X$ at the origin (resp. the blowing-up ${\rm Bl}_2(J_X)$ of $J_X$ along its 2-torsion).

\begin{proposition}
In genus 2 and characteristic 3, the "second fundamental form" map $$\widetilde{\Phi}_1 : \wH \to {\rm Bl}_2(J_X)$$ is an isomorphism.\end{proposition}
\begin{proof} It is equivalent to prove that $\Phi_1 : \mathcal{H} \to \mathcal{B}(J_{X_1})$ is an isomorphism.
 Take a point $(E,\,\alpha)$ of $\mathcal{H}$. The quotient bundle $Q_E$  is a line bundle (Lemma \ref{lemQstable} below) that is isomorphic to $$\det (F_*\O_X) \tens (\det E)^{-1}\cong \omega_{X_1}\tens (\det E)^{-1}$$ and there is an exact sequence
$$0 \to E  \xrightarrow{\alpha} F_*\O_X \to \omega_{X_1}\tens (\det E)^{-1} \to 0$$Because $\bar{\Upsilon} : \Lambda^2B \to \omega_{X_1}$ is an isomorphism in characteristic 3, the map $\widetilde{\Upsilon} : \Lambda^2(F_*\O_X) \to F_*\omega_X$ also is an isomorphism and one can identify the maps $\Phi_1(E) : \Lambda^2(E) \to F_*\omega_X$ and  $\Lambda^2(\alpha)$.

But taking second exterior powers in the short exact sequence above yields an exact sequence
$$0 \to \Lambda^2(E) \xrightarrow {\Lambda^2(\alpha)} \Lambda^2(F_*\O_X) \to E\tens \omega_{X_1}\tens \det E^{-1} \to 0$$ Because $E\tens (\det E)^{-1} \cong E^\vee$, we find that the two exact sequences correspond to each other by Serre's duality.\end{proof}

\begin{remark}\label{rkLP2Section6} This is actually the construction of \cite[Section 6]{LP2}. In particular, we find in that case that $\wH$ is smooth and that $\widetilde{\Phi}_1$ is an étale map.
\end{remark}

\subsection{A closed immersion via the quotient bundle}

In rank 2, the lemma \ref{lemsubF*xi} has the following other consequence.
\begin{lemma} \label{lemQstable}Let $(E,\,\alpha)$ be a point of $\mathcal{H}_\eta$. Then the quotient bundle fitting in the exact sequence $$0 \to E \xrightarrow{\alpha} F_*\eta \to Q_{E,\,\eta} \to 0$$ is a stable vector bundle with rank $p-2$ and degree $(p-1)(g-1)+d$, where $d=\deg \eta$.\\
In particular, the moduli property of $M_{X_1}^Q$ defines a forgetful morphism \begin{equation}\label{eqdefhM} h_\mathcal{Q} : \mathcal{H} \to \mathcal{M}_{X_1}^Q\end{equation}
\end{lemma}
\begin{proof}Assume first that $Q_{E,\,\eta}$ has torsion. Then, the kernel of the map $F_*\eta \to Q_{E,\,\eta}^{\vee\vee}$ to the torsion free part of $Q_{E,\,\eta}$ has rank 2 and slope at least $1/2$. But $$\mu_2(\eta)=(g+d-1)/p\leq (g-1)/p < 1/2$$ and this is a contradiction. If $p=3$, $Q_{E,\,\eta}$ is a line bundle and it is stable. If $p\geq 5$, let $M$ be a quotient bundle of $Q_{E,\,\eta}$ with rank $r\leq p-3$. It is also a quotient of $F_*\eta$ and its slope $\mu$ satisfies
$$\mu\geq \nu_r(\eta)\geq \nu_{p-3}(\eta) = \frac{(p+2)(g-1)+d}{p}> \frac{(p-1)(g-1)+d}{p-2} $$
for all $p\geq 5$ and $Q_{E\,\eta}$ is again stable.\end{proof}

The forgetful map $h_{\mathcal{Q}}$  maps  a $T$-point $(\mathcal{E}_T,\,\alpha_T)$ of $\mathcal{H}$ to the $T$-point $\mathcal{Q}_T$ of $\mathcal{M}_{X_1}^Q$, where $\mathcal{Q}_T$ fits in the exact sequence of vector bundles
$$0 \to \mathcal{E}_T \xrightarrow{\alpha_T} q_{X_1}^*(F_*\O_X) \to \mathcal{Q}_T \to 0$$
over $X_1\times T$.

\begin{proposition}\label{prophQclosedimm}Assume that $p>2g$. Then, the morphism $h_\mathcal{Q} : \mathcal{H} \to \mathcal{M}_{X_1}^Q$ is a closed immersion.\end{proposition}
\begin{proof} One only needs to show that, for any point $(E,\,\alpha)$ in $\mathcal{H}$, $\Hom(F_*\O_X,\,Q_{E})$ has dimension 1. Indeed, we can apply the functor $\Hom(-,\,Q_{E})$ to the exact sequence \eqref{sesdefQ} and the connecting morphism
$$\Hom(E,\,Q_{E}) \to \Ext^1(Q_{E},\,Q_{E})$$ (that agrees with the tangent map $Th^Q : T_E \mathcal{H} \to T_{Q_{E}}\mathcal{M}_{X_1}^Q$ at $E$) is injective in this case. Using adjunction and relative duality, there is an isomorphism
$$\Hom_{X_1}(F_*\O_X,\,Q_{E}) \cong \Hom_X(\omega_{X}^{p-1},\,F^*Q_{E})$$

\begin{lemma} Let $(E,\,\alpha)$ be a point of $\mathcal{H}$. With the notations of Lemma \ref{lemPsi_1=adLambda2}, define the torsion sheaf
$$k_{\Phi_1(E)} \cong  {\rm coker}({\rm ad}(\Phi_1(E)) : \det E_0 \to \omega_X)$$ over $X$. The exact sequence \eqref{sesdefQ} over $X_1$ gives rise to an exact sequence\begin{equation*}0 \to B_2 \to F^*Q_E \to k_{\Phi_1(E)} \to 0\end{equation*} over $X$ (where $B_2 = \ker (\psi_1 : B_1 \to \omega_X)$, see Subsection \ref{subsecprelimB}).\end{lemma}
 \begin{proof} Consider the exact sequence \eqref{sesF*EDelta} for some effective  divisor $\Delta$ associated to $(E,\,\alpha)$. Letting $k_\Delta$ denote the cokernel of the injective map $\O(-\Delta) \xrightarrow{s_\Delta} \O_X$ such that $\varphi_0={\rm ad}(\alpha)$ coincides with the composite $F^*E \twoheadrightarrow \O(-\Delta) \xrightarrow{s_\Delta} \O_X$, there is an exact sequence
$$0 \to \O(\Delta)\tens \det E_0 \to E_0 \xrightarrow{\varphi_0} \O_X \to k_\Delta \to 0$$ of coherent sheaves over $X$ and the Frobenius pull-back of the exact sequence \eqref{sesdefQ} induces an exact sequence
$$0 \to \O(\Delta)\tens \det E_0 \xrightarrow{\rho_1} B_1 \to F^*Q_E \to k_\Delta \to 0$$
 Since the composite $\O(\Delta)\tens \det E_0 \xrightarrow{\rho_1} B_1 \xrightarrow {\psi_1}\omega_X$ defined at \eqref{eqrho1psi1} is non zero, it is injective and $B_2$ identifies with a subbundle of ${\rm coker}(\rho_1)$, hence of $F^*Q_E$. Now, the quotient ${\rm coker}(\rho_1)/B_2$ is also the cokernel of the maps defined in \eqref{eqrho1psi1} and the lemma follows.\end{proof}

It follows from the filtration \eqref{eqfiltF*F*OX} of $B_0$ and from the previous lemma that there is an exact sequence
$$0 \to B_2/\omega_X^{p-1} \to F^*Q_E /\omega_X^{p-1} \to k_{\Phi_1(E)}\to 0$$
\begin{lemma} When $p> 2g$ and for any $(E,\,\alpha)$ in $\mathcal{H}$, the quotient sheaf $F^*Q_E /\omega_X^{p-1}$ is torsion free.\end{lemma}
\begin{proof}
 Assume that the torsion $D:={\rm Tors}(F^*Q_E /\omega_X^{p-1})$ of the quotient sheaf $F^*Q_E /\omega_X^{p-1}$ is not zero. Then the kernel of the composite map $$F^*Q_E \to F^*Q_E/\omega_X^{p-1} \to(F^*Q_E/\omega_X^{p-1})/D$$ can be written under the form $\omega_X^{p-1}\tens \O(D)$ where $\O(D)$ is a line bundle of degree $\deg D>0$ over $X$. By functoriality of relative duality, the map $F_*\O_X \to Q_E$ factors through the surjective map $F_*\O(D) \to Q_E$. But the remark \ref{rkquotF*xi} following Lemma \ref{lemsubF*xi} tells us that $$\nu_{p-2}(\O(D))=\frac{(p+1)(g-1)+\deg (D)}{p}>\frac{(p-1)(g-1)}{p-2} $$ and this is a contradiction.\end{proof}

Therefore, if $\dim \Hom(\omega_X^{p-1},\, F^*Q_E)\geq 2$, there is an injective map $\omega_X^{p-1}\oplus \omega_X^{p-1}\hookrightarrow F^*Q_E$  and the kernel of the composite
$$\omega_X^{p-1}\oplus \omega_X^{p-1}\hookrightarrow F^*Q_E \to k_{\Phi_1(E)}$$ identifies with a rank 2 subbundle of $B_2$ with slope at least $(2p-3)(g-1)$. But $B_2$ has a unique such subbundle, namely $B_{p-2}$ (see Subsection \ref{subsecprelimB} again) which is the unique non split extension of $\omega_{X}^{p-2}$ by $\omega_{X}^{p-1}$. Therefore, it cannot be a subbundle of $\omega_X^{p-1}\oplus \omega_X^{p-1}$ and this is a contradiction. \end{proof}

\subsection{Fitting ideals}

Recall
from the GIT construction of the moduli space $M_{X_1}^Q$ that there
is a Quot-scheme $Quot$, quasi-projective over $k$, that represents
the functor from the category ${\bf Sch}/k$ of schemes over $k$ to
${\bf Sets}$ defined by
$$T \mapsto \left\{k^N\tens \O_{X_1\times T} \twoheadrightarrow \mathcal{Q}_T
\tens q_{X_1}^*L_{>0}\left|\begin{array}{c}
\,\mathcal{Q}_T \text{ is locally free of rank }p-2 \text{ and,}\\
\forall t\in T(k), \mathcal{Q}_T(t) \text{ is
stable}\\
\text{and has degree }(p-1)(g-1).\end{array}\right.\right\}/ \sim $$ where $L_{>0}$ is a
suitably high degree line bundle over $X_1$. We require that there is an
isomorphism $k^N\tens \O_T \xrightarrow \sim q_{T,\,*}(\mathcal{Q}_T \tens
q_{X_1}^*L_{>0})$ that identifies a surjective map  \mbox{$k^N\tens \O_{X_1\times
S} \twoheadrightarrow \mathcal{Q}_T \tens q_{X_1}^*L_{>0}$} to the
evaluation of global sections of $\mathcal{Q}_T \tens
q_{X_1}^*L_{>0}$. This gives a natural action of $P{\rm
GL}(N)$ on $Quot$ (resp. on the universal quotient bundle
$\mathcal{Q}_{Quot}$ over $X_1\times S$) and those actions are
compatible. The coarse moduli property associates to the universal
bundle $\mathcal{Q}_{Quot}$ a morphism $\pi : Quot \to M_{X_1}^Q$
which makes $Quot$ into a locally isotrivial principal $P{\rm
GL}(N)$-bundle over the stable locus $(M_{X_1}^Q)^s$ of $M_{X_1}^Q$.

\begin{remark} If $U$ is an affine open subset in $(M_{X_1}^Q)^s$, $\pi^{-1}(U)$
is an affine open subset of $Quot$ and there is an isomorphism
$$H^0(U,\O_U) \xrightarrow \sim H^0(\pi^{-1}(U),\,\O_{Quot})^{P{\rm
GL}(N)}$$ In particular, any $P{\rm GL}(N)$-invariant sheaf of
ideals over $Quot$ descends to a sheaf of ideals over $(M_{X_1}^Q)^s$.
\end{remark}

For any $T$-point $\mathcal{Q}_T : T \to Quot$, we let
$$\mathfrak{I}^Q:={\rm Fitt}_{((p-2)^2-2)(g-1)}(R^1
q_{T,\,*}(\mathcal{H}om(F_*\O_X,\,\mathcal{Q}_T))) \subset
\O_{T}$$
denote the $((p-2)^2-2)(g-1)$-th Fitting ideal of the first direct image of $\mathcal{H}om(F_*\O_X,\,\mathcal{Q}_T)$, where $\mathcal{Q}_T$ also denotes the vector
bundle $\mathcal{Q}_T^*(\mathcal{Q}_{Quot})$. Because $X_1$ is a curve, the formation of $R^1 q_{T,\,*}$ commutes
with any base change. It is also the case for the formation of
Fitting ideals (see \cite[Corollary 20.5]{E}), and there is a canonical isomorphism $\mathfrak{I}_T
\cong \mathcal{Q}_T^*(\mathfrak{I}_{Quot})$ compatible with base
change. As a consequence of the construction, we obtain the following

\begin{lemma}The Fitting ideal $\mathfrak{I}_{Quot}$ is $P{\rm
GL}(N)$-invariant, hence descends to an sheaf of ideals
$\mathfrak{I}$ over $M_{X_1}^Q$.\end{lemma}

One checks that the sheaf $$\mathcal{F}^1_T:= q_{T,\,*}\mathcal{H}om(q_{X_1}^*(F_*\O_X\tens \omega_{X_1}^{-1}),\,\mathcal{Q}_T)$$ is locally free of rank $(g-1)(p^2-2)$ and, assuming that $\Hom(F_*\O_X,\,\mathcal{Q}_T(t))=0$ for a generic point $t$ in $T$,
one obtains a projective resolution
\begin{equation}\label{projresolFittingideal}0 \to \mathcal{F}_T^1 \to \mathcal{F}_T^0 \to R^1q_{T,\,*}Hom(q_{X_1}^*F_*\O_X,\,\mathcal{Q}_T)\to 0\end{equation}
of $R^1q_{T,\,*}Hom(p_{X_1}^*F_*\O_X,\,\mathcal{Q}_T)$, where $$\mathcal{F}_T^0:=\O_{T}^{\oplus 2(g-1)p(p-2)}$$ is a locally free of rank $2(g-1)p(p-2)$,  by taking a non zero global section of $\omega_{X_1}$. Notice that if $p>2g$, it is the case that a generic $Q$ in $(M_{X_1}^Q)^s$ does not lie in the closed subscheme $H$. Indeed, it follows from Corollary \ref{cordegDelta<g} and Lemma \ref{lemfib(wh)} that
$$\dim H \leq \dim S_{X_1}+g-1 =5(g-1)-1\leq (p-2)^2(g-1)-1$$ In particular, there is a projective resolution as above when $T=Quot$.

\begin{proposition} \label{propHdetvars} Assume that $p>2g$. Then, \\
(1) the Fitting ideal $\mathfrak{I}^Q$ is the defining ideal of the closed subscheme $\mathcal{H}$ in $ M_{X_1}^Q$;\\
(2) if $\mathcal{H}$ has the expected dimension $2(g-1)$, it is locally complete intersection (hence Cohen-Macaulay) and equidimensionnal;\\
 (3) if $(E,\,\alpha)$ is a point of $\mathcal{H}$ such that $\dim \Hom(E,\,Q_E)=2(g-1)$, then $\mathcal{H}$ has the expected dimension at $E$ and it is regular at this point.
\end{proposition}
\begin{proof} Since for any $Q$ in $M_{X_1}^Q$, a non-zero map $F_*\O_X \to Q$ is surjective (use Lemma \ref{lemsubF*xi}), it is clear that $\mathcal{H}$ and the zero locus $\mathcal{H}'$ of $\mathfrak{I}^Q$ have the same support. To prove that there is a scheme-theoretic isomorphism $\mathcal{H}\cong \mathcal{H}'$, it is enough to prove that $\mathcal{H}(T)=\mathcal{H}'(T)$ for any affine $k$-scheme $T$ (the inclusion $\mathcal{H}(T)\subset \mathcal{H}'(T)$ is clear from the definitions).

Assume that $\mathcal{Q}_T$ is a $T$-point of $M_{X_1}^Q$ that cancels the Fitting ideal $\mathfrak{I}^Q$. Using the corresponding morphism to the Quot-scheme $Quot$ (still denoted by $\mathcal{Q}_T$), one pulls-back the projective resolution of $R^1q_{Quot,\,*}Hom(q_{X_1}^*F_*\O_X,\,\mathcal{Q}_{Quot})$ and obtains an exact sequence
$$\mathcal{Q}_T^*(\mathcal{F}_{Quot}^1) \to \mathcal{Q}_T^*(\mathcal{F}_{Quot}^0) \to R^1q_{T,\,*}Hom(q_{X_1}^*F_*\O_X,\,\mathcal{Q}_{T})\to 0$$
Using Proposition \ref{prophQclosedimm} and its proof, we find that $R^1p_{T,\,*}Hom(p_{X_1}^*F_*\O_X,\,\mathcal{Q}_T)$ is locally free of rank $((p-2)^2-2)(g-1)+1$ and the kernel of $\mathcal{Q}_T^*(\mathcal{F}_{Quot}^1) \to \mathcal{Q}_T^*(\mathcal{F}_{Quot}^0)$  is locally free of rank 1. Locally on any sufficiently small open subset $U$ of $T$, one can therefore find a non-zero morphism $\beta_U : p_{X_1}^*F_*\O_X \to \mathcal{Q}_U$ which is surjective since it is so at any closed point of $U$ (use Nakayama's lemma). Since $\mathcal{H}$ is a closed subscheme of $M_{X_1}^Q$ and since $T$ is affine, the $\beta_U$ can be glued and we conclude that $\mathcal{H}'(T)\subset \mathcal{H}(T)$. This proves (1).

Assume that $\mathcal{H}$ has the expected dimension $2(g-1)$. At a point $(E,\,\alpha)$ of $\mathcal{H}$, choose an affine neighborhood $U$ of $t$ in $(M_{X_1}^Q)^s$ and a finite étale cover $\bar{U} \to U$ that trivializes the $P{\rm GL}(N)$-bundle $\pi : Quot \to
(M_{X_1}^Q)^s$. Choose a section $\bar{U} \to Quot$ and let $\mathcal{Q}_{\bar{U}}$ denote the corresponding  pull-back of the universal quotient bundle $\mathcal{Q}_{Quot}$ over $\bar{U}$. Also denote by $\mathcal{H}_{\bar{U}}$ the inverse image of $\mathcal{H}$ in $\bar{U}$. It is locally complete intersection and since $(\mathcal{M}^Q_{X_1})^s$ (hence $\bar{U}$) is smooth, $\mathcal{H}_{\bar{U}}$ is Cohen-Macaulay (see, e.g., \cite[Proposition II.8.23]{Ha}). Because it is defined as a determinantal variety, one knows (see, \emph{e.g.}, \cite[Chapter 14]{Fu}) that each irreducible component of $\mathcal{H}_{\bar{U}}$ (hence of $\mathcal{H}$) has dimension at least $2(g-1)$, hence equal to $2(g-1)$. All these notions go through étale descent (see the discussion in \cite[p. 466]{E} for locally complete intersection) and this proves (2).

 Assume that  $\dim \Hom(E,\,Q_E)=2(g-1)$. Since $\Hom(E,\,Q_E)$ is isomorphic to the tangent space of $\mathcal{H}$ at $(E,\,\alpha)$, the dimension of $\mathcal{H}$ at this point is at most $2(g-1)$. Since $\mathcal{H}$ is a determinantal variety, its dimension is at least $2(g-1)$. Now, (3) follows from standard facts on Quot-schemes of the expected dimension.
\end{proof}

\begin{remark} \label{rkH2singulier} Notice that one can do the same kind of reasoning for the image of $h_\mathcal{E}$ in $N_{X_1}$ that is characterized by $\Hom(E,\,F_*\O_X)\neq 0$. This proves in particular that if $E$ is a stable rank 2 subbundle of $F_*\O_X$ such that $F^*E\cong \O_X\oplus \O_X$, its image in $N_{X_1}$ is a singular point of the image $h_\mathcal{E}(H)$ (see \cite{Fu} again). By constrast, if $(E,\,\alpha)$ is a stable point of $\mathcal{H}$ such that $F^*E\neq \O_X \oplus \O_X$ (see Subsection \ref{subsectPrym} for the existence of such a point), applying the functor $\Hom_{X_1}(E,\,-)$ to the exact sequence \eqref{sesdefQ} yields an injection $$T_{(E,\,\alpha)}\mathcal{H}\cong \Hom_{X_1}(E,\,Q_E) \hookrightarrow H^1(X_1,\,\mathcal{E}nd(E)) \cong T_EN_{X_1}$$ and the morphism $h_\mathcal{E}$ is locally (in the neighborhood of a point in $\mathcal{H}\setminus \mathcal{H}_2$) a closed immersion.\end{remark}

\section{Distinguished loci} Theorem \ref{thmTongJ} applies and it follows from a result of Nakajima (see, {\emph e.g.}, \cite{Zh}) that $X$ is even \emph{ super-ordinary} is the sense that all of its abelian Galois étale covers are ordinary.

\subsection{The semi-stable boundaries of $\mathcal{H}$ and $\wH$}\label{subsecHss}

If $E$ a strictly semi-stable rank 2 and degree 0 vector bundle, it has a degree 0 line subbundle $M$. Now, if $E$ is a subbundle  of $F_*\O_X$, so is $M$ and the latter has order dividing $p$. Upon twisting $E$ by a suitable element of $G$, we can assume that $E$ admits $\O_{X_1}$ as a line subbundle and because of the exact sequence
\eqref{sesdefB1}, the quotient bundle $E/\O_{X_1}$ identifies with a degree 0 line subbundle of the bundle $B$ of locally exact differential forms.

If $L$ is a point of $\Theta_B$ and $L\hookrightarrow B$ is the corresponding injection, we define $E^L$ as the inverse image of $L$ in $F_*\O_X$ via the canonical surjection $F_*\O_X \to B$. There is an extension
\begin{equation}\label{sesdefE^L} 0 \to \O_{X_1} \to  E^L \to L \to 0\end{equation} (that also depends on the inclusion $L \hookrightarrow B$) and $Q_{E^L}$ is canonically isomorphic to the quotient sheaf $$Q_L:=Q_{L\hookrightarrow B}\cong B/L$$ fitting in the exact sequence
\begin{equation}\label{sesLBQL}
0 \to L \to B \to Q_L \to 0\end{equation}
Using the Quot-scheme $\mathcal{H}_B$ defined in Subsection \ref{subsecprelimB}, this defines closed immersions $\mathcal{H}_B \to \mathcal{H}$ and $\mathcal{H}_B \to \mathcal{H}/G$  and $\mathcal{H}_B$ identifies in particular with the semi-stable boundary of $\mathcal{H}/G$. Thus, the semi-stable boundary $\wH^{ss}$ of $\wH$ identifies with a degree $2^{2g}$ étale covering of $\mathcal{H}_B$.

\begin{remark}\label{rkFEsssplit} Notice that the adjoint map $F^*E^L \to \O_X$ provides with a splitting of the Frobenius pull-back of the exact sequence \eqref{sesdefE^L}. In other words, $F^*E^L\cong \O_X \oplus F^*L$.\end{remark}

\subsection{Vector bundles arising from Prym varieties}\label{subsectPrym}
We recall some basic facts about Prym varieties associated to double étale covers of $X$ and their relation with rank 2 vector bundles invariant under the action of an order 2 line bundle over $X_1$. The reader can look in \cite{Mu} for proofs and details (see also \cite{Du} where we use these technics in the same context).

\def\wX{\widetilde{X}}

Choose an order 2 line bundle $\tau$ over $X$ and consider the corresponding double
étale cover $\alpha : \wX:={\bf Spec}(\O_X\oplus \tau) \to
X$. To this cover is classically associated the homomorphisms $\alpha^* : J_X \to J_{\widetilde{X}}$ (with finite kernel isomorphic to $<\tau>$) and ${\rm Nm} : J_{\widetilde{X}} \to J_X$ defined by $M \mapsto M\tens i^*M$ (where $i : \widetilde{X} \to \widetilde{X}$ is the involution permuting the sheets of the double cover $\alpha$). The latter is surjective and its kernel is the product of a finite group isomorphic to $(\mathbb{Z}/2\mathbb{Z})$ and of a principally polarized abelian variety $P_\tau$ (the Prym variety associated to $\alpha$)  of
dimension $g-1$. The composite ${\rm Nm}\circ \alpha^*$ is $[2] : J_X \to J_X$.

For any $L$ in $J_{\widetilde{X}}$, $\alpha_*L$ is a rank 2 vector bundle with determinant ${\rm Nm}(L)\tens
\tau$, satisfying $\alpha_*L\tens \tau\cong \alpha_*L$. It is stable if  $L$ does not belong to ${\rm
im}(\alpha^*)$ and $\alpha_*\alpha^*L\cong
(\O_X\oplus \tau)\tens L$ for any $L$ in $J_X$.
Conversely, if $E$ is
any vector bundle of rank 2 and degree 0, an isomorphism $E\tens
\tau \cong E$ allows one to give $E$ a structure of $(\O_X\oplus
\tau)$-free module of rank 1, meaning that $E\cong \alpha_*L$ for some $L \in
J_{\widetilde{X}}$. In particular, there is a morphism ${\rm Nm}^{-1}(\tau) \to S_{X}$ defined by $L\mapsto \alpha_*L$ that surjects onto the $\tau$-invariant locus of $S_X$ which identifies with the intersection of $S_X$ with the (disjoint union of) two projective spaces in $|2\Theta|$ that are invariant under the action of $\tau$. It is easily checked that, restricted to an irreducible component of ${\rm Nm}^{-1}(\tau) $ (that identifies with $P_\tau$), the morphism $P_\tau \to S_X$ factors through the Kummer variety $K_{P_\tau}:=P_\tau/\{\pm 1\}$ of $P_\tau$ and that the latter is embedded in $S_X$ by this map. In particular, it intersects the strictly semi-stable locus $S_{X}^{ss}\cong K_X$ in $2^{2(g-1)}$ points that identify with the 2-torsion $P_\tau[2]=P_\tau \bigcap {\rm im}(\alpha^*)$ of $P_\tau$.

Since $p$ is odd, $F$ induces an isomorphism $J_{X_1}[2] \xrightarrow \sim J_X[2]$ and letting $\tau_1$ denote the image of $\tau$, the Frobenius inverse image $\alpha_1 : \wX_1 \to X_1$ of $\alpha : \wX \to X$ is the étale double cover associated to $\tau_1$. We choose a morphism $P_{\tau_1} \to S_{X_1}$ associated to $\tau$. Following \cite[Lemma I.11]{SGA1}, one has a cartesian diagram \begin{center}
\unitlength=0.6cm
\begin{picture}(8,4)
\put(1,2.5){$\widetilde{X}$} \put(6,2.5){$\widetilde{X}_1$}
\put(1,0){$X$}\put(6,0){$X_1$} \put(1.6,2.7){\vector(1,0){4}}
\put(1.6,0.2){\vector(1,0){4}} \put(1.2,2.3){\vector(0,-1){1.5}}
\put(6.2,2.3){\vector(0,-1){1.5}} \small \put(3.5,0.4){$F$}
\put(3.5,2.9){$F$} \put(0.5,1.4){$\alpha$} \put(6.6,1.4){$\alpha_1$}
\end{picture}\end{center} and the action of Frobenius on vector bundles arising from the Prym variety $P_{\tau_1}$ is given by multiplication by $p$ on $P_{\tau_1}$.

\begin{proposition} \label{propEprymsmooth} Let $\tau_1$ be any order 2 line bundle on $X_1$ and let $(E,\,\alpha)$ be a stable $\tau_1$-invariant point of $H$. Then, one has $\dim \Hom_{X_1}(E,\,Q_E)=2(g-1)$.\\
The scheme-theoretic locus of $\tau_1$-invariant points of $\bar{H}$ is finite reduced and contains exactly $2^{2(g-1)-1}(p^g-1)$ points.
\end{proposition}
\begin{proof} A point $(E,\,\alpha)$ in $\mathcal{H}$ such that $E\tens \tau_1\cong E$ can be written under the form $\alpha_{1,\,*}L$ for some line bundle $L$ in $J_{\widetilde{X}_1}$, and the exact sequence \eqref{sesdefQ} becomes
$$0 \to \alpha_{1,\,*}L \to F_*\O_X \to Q_{\alpha_{1,\,*}L} \to 0$$
Use adjunction and relative duality for the map $\alpha_1$ and the commutativity of the diagram preceding the statement to obtain an injective map $L \to F_*\O_{\widetilde{X}}$ and deduce the fact that $L$ has order $p$ (since $E$ is stable, $L \neq \O_{\widetilde{X}_1}$). Since $F^*\alpha_{1,\,*}L \cong \alpha_{*}\O_X \cong \O_X \oplus \tau$, a $\tau_1$-invariant point of $\bar{H}$ neither lies in the base locus of $V$ nor in $\bar{H}_2$. Since ${\rm Nm}(L)$ also has order $p$ and since $[2] : J_{X_1} \to J_{X_1}$ induces an isomorphism of $J_{X_1}[p\,]$, we can assume, upon twisting by a suitable order $p$ line bundle over $J_{X_1}$, that ${\rm Nm}(L)\cong \O_{X_1}$.

Pull the exact sequence above back to $\widetilde{X}_1$ to
obtain the exact sequence $$0 \to L^{-1}\oplus L \to F_*\O_{\widetilde{X}} \to
\alpha_1^*Q_{\alpha_{1*}M} \to 0$$ and twist the latter by $L$. Consider the
following commutative diagram of $\O_{\widetilde{X}_1}$-modules with
exact rows and columns
$$\begin{array}{ccccccc}
& 0 & & 0\\
& \downarrow & & \downarrow \\
& \O_{\widetilde{X}_1} & = & \O_{\widetilde{X}_1}\\
& \downarrow & & \downarrow \\
0 \to &\O_{\widetilde{X}_1} \oplus L^2 & \to & F_*\O_{\widetilde{X}}
& \to & \alpha_1^*Q_{\alpha_{1*}L} \tens
L & \to 0\\
& \downarrow &  & \downarrow &  & || \\
0 \to & L^{2} & \to & \widetilde{B} & \to &
\alpha_1^*Q_{\alpha_{1*}L} \tens
L & \to 0\\
& \downarrow &  & \downarrow \\
& 0 & & 0\end{array}$$ where the middle column is the exact sequence
\eqref{sesdefB1} for $\widetilde{X}_1$.
Use projection formula for $\alpha_1$ and the isomorphism $\alpha_{1,\,*}L\cong (\alpha_{1,\,*}L)^\vee$ (deriving from the assumption ${\rm Nm}(L)\cong \O_X$) to obtain an isomorphism $\alpha_{1,\,*}(\alpha_1^*Q_{\alpha_{1*}L}\tens L) \cong \mathcal{H}om_{\O_{X_1}}(\alpha_{1,\,*}L ,\,Q_{\alpha_{1*}L})$. The long exact sequence of cohomology
associated to the bottom line thus yields a surjection
$$H^1(\wX_1,\,\widetilde{B}) \twoheadrightarrow \Ext^1_{X_1}(\alpha_{1,\,*}L,\,Q_{\alpha_{1,\,*}L})$$
Because $X$ is a general curve, it is super-ordinary and its double étale cover $\wX$ is ordinary. Therefore, $\Ext^1_{X_1}(\alpha_{1,\,*}L,\,Q_{\alpha_{1,\,*}L})=0$ and the dimension of $\Hom_{X_1}(\alpha_{1,\,*}L,\,Q_{\alpha_{1,\,*}L})$ is given by Riemann-Roch. Namely, one has
 $$\dim (\Hom_{X_1}(\alpha_{1,\,*}L,\,Q_{\alpha_{1,\,*}L}))=2(p-2)\left(\frac{p-1)(g-1)}{p-2}-g+1\right)=2(g-1)$$

Now, assume that $(\mathcal{E}_T,\,\alpha_T)$ is a $T$-point of $\mathcal{H}$ for some connected $k$-scheme $T$. Assume that there is an isomorphism $\mathcal{E}_T\times q_{X_1}^*\tau_1 \xrightarrow \sim \mathcal{E}_T$ over $X_1\times T$. Then, there is a morphism $T \to P_{\tau_1}$ such that the composite $T \to P_{\tau_1} \xrightarrow{L\mapsto \alpha_{1,\,*}L} N_{X_1}$ coincide with the composite $T \to \mathcal{H} \xrightarrow{h_\mathcal{E}} N_{X_1}$. This morphism $T \to P_{\tau_1}$ factors through the reduced subscheme $P_{\tau_1}[p\,]_{\rm red}$ and it proves that the $\tau_1$-invariant locus of $\mathcal{H}$ is finite and reduced. Since $\O_{X_1}\oplus \tau_1$ cannot be imbedded in $F^*\O_X$, we find that the intersection of $\bar{H}$ and of the image of $P_{\tau_1}$ agrees with the stable part of the image of $P_{\tau_1}[2p\,]$ endowed with the reduced induced structure. Hence, it is in one-to-one correspondence with the set
$$\left(P_{\tau_1}[2p\,]\setminus P_{\tau_1}[2]\right)/\{\pm\}$$whose cardinality is $(2^{2(g-1)}p^g-2^{2(g-1)})/2=2^{2(g-1)-1}(p^g-1)$.\end{proof}

Applying Proposition \ref{propHdetvars}.(3), one has in particular
\begin{corollary}
If $p>2g$, $\wH$ is smooth at any $\tau_1$-invariant point.
\end{corollary}

\subsection{The set of Frobenius trivialized points of $\wH$}
The set of (integrable) connections on $\O_X\oplus \O_X$ is affine under $\Hom_X(\O_X\oplus \O_X,\,(\O_X\oplus \O_X)\tens \omega_X)$ and we use the trivial connection induced by the isomorphism $F^*(\O_{X_1}\oplus \O_{X_1})\cong \O_X\oplus \O_X$. The group $\Aut_X(\O_X\oplus \O_X)\cong {\rm GL}(2,\,k)$ acts on the vector space $\Hom_X(\O_X\oplus \O_X,\,(\O_X\oplus \O_X)\tens \omega_X)$ by conjugation (we will talk about \emph{transport}) and the vanishing of the $p$-curvature is invariant under transport.

Since we will perform some explicit computations, let us recall the following general technic. Being given a rank $r$ bundle $M$ with connection $\nabla$ over a curve $X$, choose an affine open subset $U$ of $X$ that trivializes $M$ and $\omega_X$, choose a
differential form $\omega_0$ (defined over $U$) that does not vanish over $U$ and provides with a trivialization of $\omega_X$. Choose a germ of derivation $\theta_0$ defined over $U$ such that
 $\theta_0(\omega_0)=1$. The restriction of $\nabla$ to $U$ is given by a $r\times r$ matrix $T$ with coefficients in
$\Gamma(U,\,\O_U)$ such that $$\nabla(f\tens e)=T(e)\tens \omega_0 + e\tens df$$ The vanishing of
$$\psi_\nabla(\theta_0)=(T+\theta_0)^p-\theta_0^p(\omega_0)T-\theta_0^p$$
suffices to prove the vanishing of the $p$-curvature. Computing $(T+\theta_0)^p$ recursively and looking for
$(T+\theta_0)^n$ under the form
$\sum_{k=0}^nT^{(n)}_k\theta_0^k$, one finds
\begin{equation}\label{eqpcurvinduction} (T+\theta_0)^{n+1}=(T+\theta_0)T_0^{(n)}+
\sum_{k=1}^n \left((\theta_0 +T)T^{(n)}_k+ T^{(n)}_{k-1}\right)
\theta_0^k + T_n^{(n)}\theta_0^{n+1} \end{equation}
By definition, one has $T_0^{(1)}=T$ and $T^{(1)}_1=1$. This implies that
$T^{(n)}_n=1$ for all $n\geq 1$ and we set $T^{(0)}_0=1$ for consistency. One checks by induction that
$$T_{n-r}^{(n)}=\left(\begin{array}{c} n \\ r
\end{array}\right)T_0^{(r)} \text{ for all }n,\,0\leq r \leq n,$$ and consequently, that $\psi_\nabla(\theta_0)=T_0^{(p)}-\theta_0^p(\omega_0)T$.

\begin{remark} \label{rkpcurvtedious} (1) The explicit computation of $T_0^{(p)}$ is  usually very tedious since $T$ and $\theta_0$ do not commute. Still, in the case of a line bundle, the expression dramatically simplifies to give
\begin{equation*}\psi_\nabla(\theta_0)
=T^p+\theta_0^{p-1}(T)-\theta_0^p(\omega_0)T\end{equation*}
(2) In case $M$ has trivial determinant, the explicit computation of the
determinant connection ensures that
$\nabla^{\rm det}$ is represented by ${\rm Tr}(T)$.\end{remark}

Let us look at a point $E$ of $\mathcal{H}$ such that $F^*E\cong \O_X\oplus \O_X$. If $E$ is stable, it is a point of $\mathcal{H}_2$. By contrast, assume that $E$ is a strictly semi-stable vector bundle that one can think of as an extension $$0 \to \O_{X_1} \to E^L \to L \to 0$$ with $L$ in $\Theta_B$ (see Subsection \ref{subsecHss}). Since $\det F^*E^L\cong F^*L \cong \O_X$, $L$ has order $p$ and as $X$ is a general curve, Theorem \ref{thmTongJ}.(3) applies and such a point $L$ is a regular point of $\Theta_B$. In particular, $\dim \Hom_{X_1}(L,\,B)=1$ and $E$ is the split vector bundle $\O_{X_1}\oplus L$ (notice, though, that $\varphi_0 : \O_X\oplus \O_X \to \O_X$ maps a local section $(f_1,\,f_2)$ to the sum $f_1+f_2$ rather than to one of the $f_i$'s, and that the second fundamental form is indeed non zero). Upon conjugating by a global automorphism of $\O_X\oplus \O_X$, there is thus a canonical $p$-integrable connection
$$\nabla^{\rm can}=\left(\begin{array}{cc} d & 0\\ 0& d+\omega_L\end{array}\right)$$ on $\O_{X}\oplus \O_{X}$, where $\omega_L$ is the unique global differential form such that $d+\omega_L : \O_X \to \omega_X$ has vanishing $p$-curvature and such that the corresponding line bundle on $X_1$ is $L$. Let $U_L$ be the affine subset of $X$ where $\omega_L$ does not vanish and let $\theta_L$ be the germ of derivation (defined on $U_L$) such that $\theta_L(\omega_L)=1$. Because the connection $d+\omega_L$ has vanishing $p$-curvature, one finds, using the $p$-curvature formula for line bundles in Remark \ref{rkpcurvtedious}, that
$$\psi_{d+\omega_L}(\theta_L)=(\theta_L+1)^p-\theta_L^p(\omega_L)-\theta_L^p=1-\theta_L^p(\omega_L)=0$$

If $x$ is any regular function defined on $U_L$, we let $\theta_L(x)$ denote the regular function $\theta_L(dx)$, where $d : \O_X \to \omega_X$ is the canonical derivation.

\begin{lemma} \label{lemdeformO+L}With the notations above, let $\omega$ be any global differential form over $X$ such that $\omega_L$ and $\omega$ are linearly independent and denote by $x$ the non constant rational function $\omega/\omega_L$. Then, none of the expressions
$$\sum_{k=1}^{p-1}\theta_L^k(x)\text{
and }\sum_{k=1}^{p-1} \left(\begin{array}{c} p-1\\ k\end{array}\right) \theta_L^k(x)$$ vanishes.
\end{lemma}
\begin{proof} By assumption, $\dim \Hom(L,\,B)=1$ for all line bundle $L$ of order $p$ on $X_1$ and the map $\Ext^1(L,\,\O_{X_1}) \xrightarrow{H^1(F^*)} \Ext^1(\O_X,\,\O_X)$ is injective. Therefore, if a $E$ is S-equivalent to $\O_{X_1}\oplus L$ and satisfies $F^*E\cong \O_X\oplus \O_X$, it is necessarily isomorphic to $\O_{X_1}\oplus L$. In terms of connections, it means that both connections
$$\nabla=\left(\begin{array}{cc} d & \omega\\ 0 & d+\omega_L  \end{array}\right) \text{ and } \nabla'=\left(\begin{array}{cc} d+\omega_L & \omega\\ 0 & d \end{array}\right)$$ on $\O_X\oplus \O_X$
have non-zero $p$-curvature since none of them is transport equivalent to $\nabla{\rm can}$. Because $$\left(\begin{array}{cc}0 & x\\ 0& 1\end{array}\right)\left(\begin{array}{cc}0 & f \\ 0& 1\end{array}\right)=\left(\begin{array}{cc}0 & x\\ 0& 1\end{array}\right) \left(\text{ resp. }\left(\begin{array}{cc}1 & x\\ 0& 0\end{array}\right)\left(\begin{array}{cc}1 & f\\ 0& 0\end{array}\right)=\left(\begin{array}{cc}1 & f\\ 0& 0\end{array}\right)  \right)$$for any rational function $f$ on $X$,
 we compute by induction (using the expression \eqref{eqpcurvinduction}) that
$$\psi_{\nabla}(\theta_L) = \left(\begin{array}{cc} 0 & \sum_{k=1}^{p-1} \theta_L^k(x)\\ 0 & 0  \end{array}\right) \left(\text{ resp. }\psi_{\nabla'}(\theta_L) = \left(\begin{array}{cc} 0 & \sum_{k=1}^{p-1} \left(\begin{array}{c} p-1\\ k\end{array}\right) \theta_L^k(x)\\ 0 & 0  \end{array}\right)\right) $$
and conclude the proof.\end{proof}

\begin{proposition} \label{propH2proper} The set $\mathcal{H}_2$ (resp. $\wH_2$) is proper.
\end{proposition}\begin{proof} Assume that the contrary holds. Then there is a point of type $\O_{X_1}\oplus L$ as above that can be continuously deformed into a point of $\mathcal{H}_2$. In particular, letting $X_\varepsilon$ denote $X\times_k \Spec k[\varepsilon]/\varepsilon^2$, the connection $\nabla^{\rm can}$ above  has a non-trivial linear infinitesimal deformation $\nabla_\varepsilon :=\nabla^{\rm can}+\varepsilon T$ over $X_\varepsilon$ (with $T$ in $\Hom_X(\O_X\oplus \O_X,\,(\O_X\oplus \O_X)\tens \omega_X)\cong H^0(\omega_X)^{\oplus 4})$) whose $p$-curvature vanishes and whose second fundamental form is non-zero (and remains so after conjugation by any element of ${\rm Aut}_{X_\varepsilon}(\O_{X_\varepsilon}\oplus \O_{X_\varepsilon}) \cong {\rm GL}(2,\,k[\varepsilon]/\varepsilon^2)$).

Because $X$ is ordinary, the kernel of $V_J : J_{X_1}\to J_X$ is reduced and such a deformation necessarily induces a trivial deformation of the determinant connection. It means that $T$ is a traceless morphism (see Remark \ref{rkpcurvtedious}.(2) above) that we write under the form
$$T=\left(\begin{array}{cc} \omega_{11} & \omega_{12}\\ \omega_{21}& -\omega_{11}\end{array}\right)=\left(\begin{array}{cc} 2\omega_{11} & \omega_{12}\\ \omega_{21}& 0\end{array}\right)-\omega_{11} \left(\begin{array}{cc} 1 & 0\\ 0& 1\end{array}\right)=T'-\omega_{11} I$$
Write $f^{(i)}$ for the function $\omega_i/\omega_L$, $i=11,\, 12$ or $21$. It is regular over $U_L$. Let $\nabla'_\varepsilon$ denote $\nabla_\varepsilon+\varepsilon \omega_{11}I=\nabla^{\rm can}+\varepsilon T'$. It follows from \cite[Corollary 3.6(iii)]{Oss} that one has
\begin{equation}\label{eqOssCor3.6iii}\psi_{\nabla_\varepsilon}(\theta_L)=
\psi_{\nabla'_\varepsilon}(\theta_L)+\varepsilon ((f^{(11)})^p+
\theta_L^{p-1}(f^{(11)})-(f^{(11)}))I\end{equation}
where we use $\theta_L^p(\omega_L)=1$ for the last term.

 Introduce the matrices $J=\left(\begin{array}{cc} 0 & 0\\ 0& 1\end{array}\right)=J^2$ and $R_0^{(n)}$, where $R_0^{(1)}=R=\left(\begin{array}{cc} 2f^{(11)} & f^{(12)}\\ f^{(21)} & 0\end{array}\right)$ and where $R_0^{(n+1)}=(\theta_L+J)R_0^{(n)} + RJ$ for all $n\geq 1$. One checks that
 \begin{eqnarray*}
(\theta_L+J+\varepsilon R)(J+\varepsilon R_0^{(n)})& =& J+\varepsilon R_0^{(n+1)}
\end{eqnarray*}
and writing $R^{(n)}_0$ under the form $\left(\begin{array}{cc} 2f^{(11)}_n  & f^{(12)}_n\\ f^{(21)}_n & 0\end{array}\right)$, we check by induction on $n\geq 1$ that
\begin{eqnarray*}
f^{(11)}_n & =& \theta_L^{n-1}(f^{(11)})\\
f^{(12)}_n & = & \sum_{k=0}^{n-1}\theta_L^k(f^{(12)})\\
f^{(21)}_n & =& \sum_{k=0}^{n-1} \left(\begin{array}{c} n-1\\ k\end{array}\right) \theta_L^k(f^{(21)})\end{eqnarray*}
Using \eqref{eqOssCor3.6iii} and the fact that $\psi_{\nabla_\varepsilon'}(\theta_L)=\varepsilon (R_0^{(p)}- R)$,
the vanishing of the $p$-curvature $\psi_\nabla(\theta_L)$ yields that $f^{(11)}$ and that both
$$\sum_{k=1}^{p-1}\theta_L^k(f^{(12)}) \text{ and }\sum_{k=1}^{p-1} \left(\begin{array}{c} p-1\\ k\end{array}\right) \theta_L^k(f^{(21)})$$
are 0. Using the lemma above, both $\omega_{12}$ and $\omega_{21}$ are colinear with $\omega_L$ and  $\nabla_\varepsilon$ is conjugate to the trivial deformation by some constant automorphism of $\O_{X_\varepsilon}\oplus \O_{X_\varepsilon}$. This is a contradiction and it proves the proposition.
\end{proof}

\section{The genus 2 case : Proof of the theorem \ref{thmbarHintegral}}

In this section, we will assume that $p\geq 5$ (recall that \cite[Section 6]{LP2} gives a complete description of $\bar{H}$ in characteristic 3).

In genus 2, $\bar{H}$ has dimension 2 since it is also defined as the residual divisor of $K_{X_1}$ in $V_S^{-1}(K_X)$ (notice that since $X$ is ordinary, the multiplicity of $K_{X_1}$ in $V_S^{-1}(K_X)$ is 1). It immediately follows from the proposition \ref{propH2proper} that the set $\bar{H}_2$ is finite (it is proper in the complement of an hypersurface in $\P^3$). Also, the base locus of Verschiebung is finite and, as a consequence, the set of 1-dimensional fibres of the morphism $\widetilde{h} : \wH \to \bar{H}$ is finite (see Lemma \ref{lemfib(wh)}). Therefore, $\wH$ also has dimension 2 and it is locally complete intersection and equidimensional (Proposition \ref{propHdetvars}).

\subsection{The semi-stable boundary of $\wH$ in genus 2} A generic point $(E,\,\alpha)$ in $\mathcal{H}$ gives rise to an exact sequence
$$0 \to F^*\det E\to F^*E \to \O_X \to 0$$
Push it forward by Frobenius to obtain an exact sequence
$$0 \to \det E \tens F_*\O_X \to E\tens F_*\O_X \to F_*\O_X \to 0$$ and notice that the composite $E \xrightarrow{1_E\tens F^*} E\tens F_*\O_X \to F_*\O_X$, agrees with $\alpha$. In  other words, there is a commutative diagram with exact rows and columns
\begin{equation}\label{diagF^*Egen}
\begin{array}{ccccccc}
& & & 0 & & 0\\
& & & \downarrow & & \downarrow \\
& & & \det E \tens F_*\O_X & =& \det E \tens F_*\O_X \\
& & & \downarrow & & \downarrow \\
0 \to & E & \to & E\tens F_*\O_X & \to & E\tens B & \to 0\\
& || & & \downarrow & & \downarrow & \\
0 \to & E & \to & F_*\O_X & \to & Q_E & \to 0\\
& & & \downarrow & & \downarrow \\
& & & 0 & & 0
\end{array}\end{equation}
It follows from Remark \ref{rkecupphi1=0} that if $E$ is a stable bundle in $\bar{H}$ with $(\det E)^2\neq \O_{X_1}$, then $F^*E$ is split (see also \cite[Remark 6.2]{LP1} for another argument involving the hyperelliptic involution). In particular, $\Hom(E,\,\det E\tens F_*\O_X)$ is non zero for such a point (hence for all by semi-continuity). Applying the functor $\Hom(E,\,-)$ to the rightmost column in the diagram above, one finds that  $H^0(\mathcal{E}nd(E)\tens B)$ is non zero
for all point of $\wH$.

\begin{lemma} \label{lemdefSigmaB+ssslocus} The set $\Sigma_B$ defined as
$$\Sigma_B=\left\{\begin{array}{c}\text{rank }2\text{ semi-stable vector bundles } E \text{ with
trivial}\\
\text{ determinant over }X_1\text{ such that }H^0(\mathcal{E}nd(E)\tens B)\neq 0\end{array}\right\} $$
is
compatible with S-equivalence. It defines a closed subset of $S_{X_1}$ whose intersection with $K_{X_1}$ agrees with the support of the divisor on $K_{X_1}$ induced  by the totally symmetric divisor $[2]^{-1}(\Theta_B)$ on $J_{X_1}$.\\
As sets, there is an inclusion $\bar{H} \subset \Sigma_B$.
\end{lemma}
 \begin{proof} Take a strictly semi-stable vector bundle $E$ with
trivial determinant written as an extension $0 \to M^{-1} \to E
\to M \to 0$ where $M$ is a degree 0 line bundle over $X_1$. We
tensor this exact sequence by $E\tens B$, $M\tens B$ and $M^{-1}\tens B$
successively and take the associated exact sequences of cohomology
spaces. Because $X$ is ordinary, $H^0(B)=H^1(B)=0$ and one finds an
exact sequence
\begin{equation*}0 \to H^0(M^{-2}\tens B) \to
H^0(\mathcal{E}nd(E)\tens B) \to H^0(M^2\tens B) \to H^1(M^{-2}\tens
B)\end{equation*} On the other hand, $\mathcal{E}nd(M\oplus M^{-1})\cong
\O_{X_1}^{\oplus 2}\oplus M^{-2}\oplus M^2$. Because the divisor
$\Theta_B$ is symmetric, we find that $H^0(\mathcal{E}nd(E)\tens B)\neq 0$ is
and only if $M^2$ is a point on $\Theta_B$, or equivalently, if
$H^0(\mathcal{E}nd(M\oplus M^{-1})\tens B)\neq 0$. Because $\Theta_B$ is totally symmetric in the sense of Mumford, it descends to a divisor on $K_{X_1}$ with support the set $$\{[M\oplus M^{-1}] \in K_{X_1}| M^2\in \Theta_B\}$$
The last assertion just restates the conclusion of the discussion preceding the statement.\end{proof}

It follows from the Theorem \ref{thmTongJ}.(1) that $\Theta_B$ is smooth and that the same holds for $\mathcal{H}_B \cong \Theta_B$. The inverse image of $\Theta_B$  via the multiplication map $[2] : J_{X_1} \to J_{X_1}$ is also smooth (and connected). Since the ramification locus of the Kummer map $J_{X_1} \to K_{X_1}$ is $J_{X_1}[2]$ and since $\Theta_B$ does not go through the origin, we find that the divisor in $K_{X_1}$ induced by $[2]^{-1}\Theta_B$ is connected and smooth, hence irreducible, with support $\Sigma_B\bigcap K_{X_1}$.

\begin{corollary} The semi-stable boundary $\bar{H}\bigcap K_{X_1}$ agrees with $\Sigma_B\bigcap K_{X_1}$. In particular, it is irreducible\end{corollary}
\begin{proof} We already have the inclusion of sets $\bar{H}\bigcap K_{X_1}\subset \Sigma_B\bigcap K_{X_1}$. Conversely, take a line bundle $M$ over $X_1$ and assume that $M^2$ lies in $\Theta_B$. Then, the inverse image of $M^2$ via the canonical surjection $F_*\O_X \to B$ provides with an extension
$0 \to \O_{X_1} \to E^{M^2} \to M^2 \to 0$ and twisting by $M^{-1}$, we obtain an exact sequence
$$0 \to M^{-1} \to E^{M^2} \tens M^{-1} \to M \to 0$$
that is non split as soon as $M^{2p}\neq \O_{X_1}$. By definition, this is a point of $\wH$ and $[M\oplus M^{-1}]$ lies in $\bar{H}$.\end{proof}

\begin{proof}[Proof of the Theorem \ref{thmbarHintegral}]
Any irreducible component of $\wH$ is mapped onto a (irreducible) 2-dimensional  closed subscheme in $\bar{H}$. In particular, it necessarily meets the semi-stable boundary $K_{X_1}$ along the 1-dimensional irreducible closed subscheme $\bar{H}\bigcap K_{X_1}$. The same holds in $\wH$. Namely, any irreducible component of $\wH$ meets a irreducible component of the (smooth) semi-stable boundary $\wH^{ss}$ and the latter identifies with a degree
16 étale cover of $\mathcal{H}_B \cong \Theta_B$. The composite $\mathcal{H}_B\hookrightarrow \mathcal{H} \xrightarrow{h_\mathcal{E}} N_{X_1} \xrightarrow \det J_{X_1}$ identifies with the composite $\mathcal{H}_B \xrightarrow \sim \Theta_B \subset J_{X_1}$ and it follows from the cartesian diagram \eqref{diagNXSXJX} and the connectedness of the inverse image $[2]^{-1}\Theta_B$ that $\wH^{ss}$ is connected, hence irreducible. As a consequence, $\wH$ is irreducible and so is $\bar{H}$.

Also, if the singular locus of $\bar{H}$ had dimension 1, it would necessarily intersect $\bar{H}\bigcap K_{X_1}$ and this is a contradiction. Since $\bar{H}$ is complete intersection in $S_{X_1}$, it is normal. The fact that it has degree $2(p-1)$ is a straightforward consequence of the proposition \ref{propEprymsmooth} and the equality of divisors follows since $V^{-1}(K_X)$ has degree $4p$.\end{proof}

\subsection{The hyperelliptic involution} In the previous subsection, we used the Remark \ref{rkecupphi1=0} to derive the fact that, for a generic point $(E,\,\xi,\,\alpha)$ of $\wH$, there is also a point $(E,\,\xi^{-1},\,\alpha')$ in $\wH$. More generally, recall that the hyperelliptic involution $\imath$ of $X_1$ changes a degree 0 line bundle over $X_1$ to its inverse. For any $T$-point $\mathcal{E}_T \xrightarrow{\alpha_T} F_*\xi_T$ of $\wH$, we define the $T$-point ($\imath^*(\mathcal{E}_T),\,\xi_T^{-1},\,\imath(\alpha_T))$ of $\wH$ as the composite
$$\imath(\alpha_T) : \imath^*\mathcal{E}_T \xrightarrow {\imath^*\alpha_T} \imath^*(F_*\xi_T)\xrightarrow \sim F_*\xi^{-1}_T$$
This defines an involution of $\wH$ (still denoted by $\imath$) whose fixed locus consists in the closure of the set of points $(E,\,\tau,\,\alpha)$ such that $\tau^2\cong \O_X$ and such that $\alpha$ induces an exact sequence
$$0 \to \tau \to F^*E \to \tau \to 0$$

We let $[\wH]$ denote the quotient of $\wH$ under this involution and we notice that $\widetilde{h} : \wH \to \bar{H}$ (resp. $\widetilde{\Phi}_1 : \wH \to {\rm Bl}_2(J_X)$) induces a morphism $[\widetilde{h}] : [\wH] \to \bar{H}$ (resp. $[\widetilde{\Phi}_1] : [\wH] \to {\rm Bl}_2(K_X)$, where ${\rm Bl}_2(K_X)$ denotes the smooth K3-surface obtained as the blowing-up of $K_X$ along its singular locus). It immediately follows from the discussion above and from Lemma \ref{lemfib(wh)} that the map $[\widetilde{h}] : [\wH] \to \bar{H}$ identifies \emph{set-theoretically} with the blowing-up of $\bar{H}$ along the union of the finite sets consisting in the base locus of $V$ on the one hand and in the points of $\bar{H}_2$ on the other hand. The exceptional line above a point in the former set canonically identifies with $X_1/\imath\cong |\omega_{X_1}|\cong \P^1$ whereas the exceptional line above a point in the latter set is isomorphic to $\P \Hom(E,\,F_*\O_X)$ since it lies in the fixed locus of $\imath$. Since $\bar{H}$ is normal, Zariski's main theorem readily gives the following

\begin{proposition} The morphism $[\widetilde{h}] : [\wH] \to \bar{H}$ identifies with the blowing-up of $\bar{H}$ along the union the finite sets consisting in the base locus of $V$ on the one hand and in the points of $\bar{H}_2$ on the other hand.\end{proposition}

\end{document}